\documentclass[11pt]{paper}
\usepackage[a4paper, margin=25mm]{geometry}
\usepackage[english]{babel}

\usepackage[utf8x]{inputenc}
\usepackage{amsmath, amsthm, amssymb}
\usepackage{mathtools}
\usepackage{enumitem}
\usepackage[nobysame,initials]{amsrefs}
\usepackage{nameref,hyperref,cleveref,cite}
\usepackage{dsfont,bbm}
\usepackage{graphicx}
\usepackage{xcolor}
\usepackage{textcomp}
\usepackage{comment}

\usepackage{wrapfig}
\usepackage{pdfpages}
\usepackage{comment}

\usepackage{array}
\newcolumntype{C}[1]{>{\centering\arraybackslash}p{#1}}

\usepackage{multirow}

\AtBeginDocument{%
  \def\MR#1{}%
}

\theoremstyle{plain}

\newtheorem{theorem}{Theorem}[section]
\newtheorem{lemma}[theorem]{Lemma}

\newtheorem*{corollary*}{Corollary}

\theoremstyle{definition}

\Crefname{claim}{Claim}{Claims}
\crefname{appendix}{Appendix}{Appendices}
\Crefname{appendix}{Appendix}{Appendices}

\newcommand{\EE}{\mathbb{E}}
\newcommand{\PP}{\mathbb{P}}
\newcommand{\VV}{\mathbb{V}}
\newcommand{\Po}{\textnormal{Po}}
\newcommand{\indi}{\mathbf 1}

\DeclareMathOperator{\dTV}{d_{TV}}

\DeclareMathOperator{\Vol}{Vol}
\DeclareMathOperator{\diam}{diam}

\DeclareMathOperator{\inter}{int}

\newcommand{\conv}{{\rm conv }\,}
\newcommand{\cl}{{\rm cl }\,}

\renewcommand{\l}{\left}
\renewcommand{\r}{\right}

\newcommand{\dd}{{\rm d}}
\newcommand{\eps}{\varepsilon}

\newcommand{\caseif}{\text{if }}

\newcommand{\R}{\mathbb{R}}
\newcommand{\N}{\mathbb{N}}
\newcommand{\Z}{\mathbb{Z}}

\renewcommand{\c}{\mathcal}
\renewcommand{\bf}{\mathbf}

\author{Kinga Nagy\footnote{Osnabrück University, Osnabrück, Germany. \newline Email: kinga.nagy@uni-osnabrueck.de}}
\title{On the Maximal Dimension of Random Simplicial Complexes}
\begin{document}

\maketitle

\begin{abstract}
    The dimension of random simplicial complexes (defined as the maximal dimension among all faces) is a natural extreme value associated with the complex, and is closely related to other functionals defined by a maximum, such as the clique number of geometric graphs or scan statistics. We extend existing results in the binomial point process case to the Poisson setting in sparse graphs, give new ones about expectations and large deviation principles in all regimes, as well as give a first precise distribution result in the dense case.  
\end{abstract}



\section{Introduction}

Ever since the introduction of the Gilbert graph in 1961 \cite{Gilbert}, stochastic models of geometric graphs have been heavily investigated. In the Gilbert graph, or random geometric graph, the vertices are given by a random point configuration in Euclidean space, and a pair of distinct vertices is connected by an edge whenever their distance is at most some threshold $r>0$; the point configuration is most often given by a binomial- or Poisson point process. Some areas of interest are then graph-theoretical or metric properties of this graph: expected number of edges or total edge length, chromatic number or clique number, connectedness. 
Much is known about each of them; the seminal book of Penrose \cite{PenroseBook} contains an introduction to the topic, as well as a large collection of results. 
The assertions are generally about asymptotic behaviour, as the number of (expected) points tends to infinity while the threshold $r$ tends to zero.

%Many developments have been made since, REF  \cites{B15,BP16,BR17,RS13,RST17,RRvW24}
Subgraphs counts, volume power functionals such as total edge length, and more generally, related $U$-statistics, have been heavily investigated in the last decades, yielding results regarding moments, concentration of the probability distribution, limit theorems and multivariate limit theorems; a non-exhaustive list of such results is \cite{AR20},\cite{B16},\cite{BP16},\cite{BR17},\cite{LRP13},\cite{RRvW24},\cite{RS13},\cite{RST17}.
As a contrast to functionals that arise as sums over certain tuples, one could consider functionals related to extremes, such as the maximal degree or clique number (i.e. size of the largest clique) of the graph. Results of this type, compared to ones about $U$-statistics, are few and far between. The maximal degree, clique- and chromatic numbers of the graph are discussed in \cite[Section~6]{PenroseBook}, \cite{M08},\cite{MM11}. The existing literature primarily deals with sparse graphs that arise from models with small distance parameter $r$, such as in \cite{M99},\cite{M08},\cite{P02}; the details of this will be discussed later. 
From a slightly different perspective, but more generally, Lao and Mayer \cite{LM08} introduced the phrase \textit{$U$-max statistic} for the extremal analog of a $U$-statistic: given a score function defined on $k$-tuples, consider its maximum over all $k$-tuples of the process. 
Their results were significantly extended in \cite{ST12}, where order statistics of score functions were considered.

As an extension of the Gilbert graph, random simplicial complexes have also been considered recently: a natural example of this is the clique complex of the Gilbert graph, that is, the simplicial complex given by its complete subgraphs, also known as the Vietoris--Rips complex. A functional of interest in this context is the dimension of the complex: this is defined as the maximal dimension that occurs among its faces. It appears that results for this functional have not been stated explicitly so far; however, in certain constructions, the dimension coincides with other well-known values. In particular, for the clique complex, the dimension is exactly one less than the clique number (the difference being simply due to how the dimension of simplices is defined). In another construction we consider, the \v Cech complex, the dimension is one less than a scan statistic; these definitions are given later on. 
In this paper, we aim to prove new results about expectations, large deviation principles, and a first result about the precise distribution for dense settings.

Let us now introduce the model precisely; a more detailed introduction to some of the concepts is given in Section~\ref{section:prelim}.
For a finite point set $X\subset \R^d$ and parameter $r>0$, the \textit{geometric graph} $\c G(X,r)$ is the graph whose vertices are the points from $X$, and a pair of distinct vertices $x,y\in X$ are connected by an edge if and only if $\Vert x-y\Vert \leq r$, with $\Vert \cdot \Vert$ denoting the Euclidean norm.
Now, fix a compact convex set $W\subset \R^d$ with unit volume. Let $\eta_t$ be a homogeneous Poisson point process on $W$ with intensity $t$, and $r_t>0$ a distance parameter depending on $t$. Then we call $\c G_t\coloneq \c G(\eta_t,r_t)$ a \textit{random geometric graph} or \textit{Gilbert-graph}. We will assume that $(r_t)_{t>0}$ is monotone decreasing, and $r_t\to 0$ as $t\to \infty$; the rate of this convergence, and in particular, the behaviour of the quantity $\rho_t\coloneq tr_t^d$ controls the denseness of the graph.
Note that while in this paper we focus on the Poisson case, the case of the binomial point process -- where the set of vertices is given by an independently chosen sample of uniform points in $W$ -- is just as well-studied.
We also mention that while we work with the Euclidean distance, most of the proofs are easily adapted to general metrics. In particular, the maximum  (or $l_\infty$) norm case has been considered instead many times, as explicit computations are often easier to handle.
%$\c G_t$ as follows: the vertices of the graph are given by $\eta_t$, and a distinct pair of points $x,y\in \eta_t$ forms an edge if and only if $\Vert x-y\Vert \leq r_t$, where $r_t>0$ is a fixed threshold depending on $t$. % i.e. $\c G_t = \l( \eta_t, \{(x,y)\in \eta_{\neq}^2 \colon \Vert x-y\Vert \leq r_t\}\r)$.

Denote by $\eta_{t,\neq}^k$ the set of all $k$-tuples of distinct points of $\eta_t$. 
The \textit{Vietoris--Rips complex} $\c V = \c V(\eta_t, r_t)$ is defined as the clique complex of $\c G_t$, that is, the tuple $\{x_0,\ldots, x_{k}\} \in \eta_{t,\neq}^{k+1}$ is in $\c V$ whenever the points form a complete subgraph of $\c G_t$. 
Equivalently, $\{x_0,\ldots, x_{k}\}$ is a $k$-dimensional face (or \textit{$k$-face}) of $\c V$ if and only if $\Vert x_i-x_j\Vert\leq r_t$ for every $i,j=0,\ldots, k$. On the other hand, $\{x_0,\ldots, x_k\}$ is a simplex in the \textit{\v Cech complex} $\c C = \c C(\eta_t, r_t)$ if and only if the points are contained in a ball of radius $r_t/2$, i.e. there is $y\in \R^d$ such that $\Vert x_i-y\Vert \leq r_t/2$ for all $i=0,\ldots, k$. 
This can also be formulated using a slightly more geometric language: $\c V$ contains $k$-tuples with \textit{diameter} at most $r_t$, and $\c C$ with \textit{circumradius} at most $r_t/2$. Both of these values can be interpreted as a measure of closeness of the point set $\{x_0,\ldots, x_k\}$. %As such, one could also think of $\c V$ and $\c C$ as sets of clusters, which are defined as 'close' sets of points. 
Note that every simplex of the \v Cech complex is also in the Vietoris--Rips complex; the converse generally does not hold. (However, it does in the $l_\infty$ case.)

The set of $k$-faces of the Vietoris--Rips complex and the \v Cech complex is denoted by $\c F_k^{\c V}$ and $\c F_k^{\c C}$, respectively, and their number by $f_k^{\c V}$ and $f_k^{\c C}$. 
Though this is omitted from the notation, they do depend on $t$ and $r_t$.
The \textit{dimension} of a simplicial complex is defined as the highest  dimension among all faces. In particular, we write $D_t^{\c V} \coloneq \max\{k\colon f_k^{\c V}>0\}$ and $D_t^{\c C} \coloneq \max\{k\colon f_k^{\c C}>0\}$ for the dimension of the Vietoris--Rips and \v Cech complexes, respectively. 
Importantly, as mentioned before, $D_t^{\c V}+1$ is the clique number of $\c G_t$, and $D_t^{\c C}+1$ is the maximal number of points in some translate of $(r_t/2)B^d$, called a \textit{scan statistic}: formally, for a Borel set $A\subset \R^d$, we call $\sup_{x\in \R^d} \eta_t(x+A)$ a scan statistic with scanning set $A$. We then have, by definition of the \v Cech complex, that $D_t^{\c C}+1$ is a scan statistic with scanning set $(r_t/2)B^d$. 
We note that unless stated otherwise, omitting the subscript %- or superscript 
$\c C$ or $\c V$ means that a statement holds for both constructions; we will emphasise the possible differences whenever necessary.

In general, we can think of the complexes $\c V$ and $\c C$ as sets of clusters, where each cluster is a set of points that are close to each other in some geometric sense. Using this viewpoint, the dimension can be considered (one less than) the size of the largest cluster in $\eta_t$.

Throughout the paper, we will primarily be using the language of simplicial complexes; however, in many cases, it is useful to think of the connection to the underlying graph, the clique number, as well as scan statistics and clusters, especially when considering the existing literature. 
%For either complex, it holds by definition that $D_t\geq k$ is equivalent to $f_k\geq 1$. 

%$g=o(h)$ if $\lim_{x \to \infty} g(x)/h(x) = 0$; $g=O(h)$ if $ \limsup_{x\to \infty} g(x)/h(x) <\infty$, and $g = \Omega(h)$ if $\liminf_{x\to \infty} g(x)/h(x)>0$. Furthermore, if $g=O(h)$ and $h=O(g)$, we say that $g=\Theta(h)$, and $g/h\to 1$ is denoted by $\sim$.}

Our main results are the following.

\begin{theorem}\label{thm:main}
Let $\eta_t$ be a homogeneous Poisson point process on the unit volume compact convex set $W\subset \R^d$ with intensity $t>0$. Choose $r_t>0$ such that $r_t\to 0$ as $t \to \infty$, and let $\rho_t\coloneq tr_t^d$. Let $D_t$ denote the dimension of the Vietoris--Rips or \v Cech complexes built on $\eta_t$ with distance parameter $r_t$. Then the following holds.
\begin{enumerate}[label=\roman*)]
    \setlength\itemsep{0pt}  %\vspace{-5pt}
    \item If $\rho_t = O(t^{-\alpha})$ for some $\alpha>0$, then there exists a unique $k\in \N$ such that $\EE f_{k-1}\to \infty$ and $\EE f_k \to \lambda\in[0,\infty)$ as $t\to \infty$. Letting $Z$ be a random variable with $\PP(Z=k-1)=e^{-\lambda}=1-\PP(Z=k)$, it holds for the total variation distance $\dTV$ that
    \[\dTV(D_t,Z) = O\big( (t\rho_t^{k-1})^{-1}\big)+O\big( t\rho_t^{k+1} \big) + O(\EE f_k-\lambda) = o(1),\]
    and in particular, $\PP(D_t\in\{k-1,k\})\to 1$. Furthermore, for any $m\in \N$, 
    \[\big| \EE D_t^m - \EE Z^m \big| = O\big( (t\rho_t^{k-1})^{-1}\big)+O\big( t\rho_t^{k+1} \big) + O(\EE f_k-\lambda) = o(1).\]
    %where clearly $\EE Z^m = (k-1)^m e^{-\lambda} + k^m (1-e^{-\lambda})$.
    %\[\big| \EE D_t^m - \EE Z^m \big| = \big|\EE D_t^m - ((k-1)^m e^{-\lambda} + k^m (1-e^{-\lambda}) )\big|\[
    % \item If $\rho_t = O(t^{-\alpha})$ for some $\alpha>0$, then there exists a unique $k\in \N$ such that $\EE f_{k-1}\to \infty$ and $\EE f_k \to \lambda\in[0,\infty)$ as $t\to \infty$, and it holds that 
    % \[\PP(D_t\in\{k-1,k\})\to 1 \text{ as } t\to\infty,\]
    % specifically $\PP(D_t=k-1)\to e^{-\lambda}$ and $\PP(D_t=k)\to 1-e^{-\lambda}$.
    % Furthermore, %for any real-valued function $g=g(x)$ for which $g=O(x^n)$ with some $n\in \R$, we have \[\lim_{t \to \infty} \EE g(D_t) = g(k-1)e^{-\lambda} + g(k)(1-e^{-\lambda}).\]
    % %As a special case, 
    % %for the moments of $D_t$ we obtain 
    % \[\EE D_t^m \to (k-1)^m e^{-\lambda} + k^m (1-e^{-\lambda}) \text{ as } t\to\infty.\] %as $t \to \infty$.
    \item If $\rho_t/t^{-\alpha}\to \infty$ for all $\alpha>0$, but $\rho_t / \ln (t) \to 0$, then there exists $k=k_t\in \N$ such that
    \[\PP(D_t\in\{k_t-1,k_t\})\to 1 \text{ and } \EE D_t^m \sim k_t^m \text{ as } t\to\infty,\]
    where \[k_t \sim \frac{\ln (t) }{\ln(\ln (t) / \rho_t)} \to \infty \text{ as } t\to \infty.\] % such that $\PP(D_t\in \{k_t-1,k_t\})\to 1$. Further, $\EE D_t^m \sim ...$
    \item If $\rho_t/ \ln (t)  \to B\in (0,\infty)$, then 
    \[D_t \sim \beta \cdot \frac{\kappa_d}{2^d} \rho_t \text{ in probability, and } \EE D_t^m \sim \l(\beta \cdot \frac{\kappa_d}{2^d} \rho_t\r)^m \text{ as } t\to \infty\]
    where $\kappa_d$ denotes the volume of the $d$-dimensional unit ball, and $\beta>1$ is the unique number satisfying $\beta \ln \beta - \beta + 1 = (B\cdot \kappa_d 2^{-d})^{-1}$. 
    \item If $\rho_t / \ln (t) \to \infty$, then 
    \[D_t \sim \frac{\kappa_d}{2^d}\rho_t \text{ in probability, and } \EE D_t^m \sim \l(\frac{\kappa_d}{2^d}\rho_t\r)^m \text{ as } t\to \infty,\] where $\kappa_d$ denotes the volume of the $d$-dimensional unit ball.
\end{enumerate}
\end{theorem}

The symbol $\sim $ means that the ratio of the quantities tends to $1$.
A fact of particular note is that the expression $ \rho_t\cdot \kappa_d 2^{-d}$ appearing in the last two regimes is exactly the expected number of points in an arbitrary ball of radius $r_t/2$ within $W$; clearly, the number of points in any such fixed ball is a lower bound for the dimension (in either complex). It is quite natural to expect that $D_t/\rho_t\to \infty$ would hold with high probability; however, somewhat surprisingly, in the denser regimes the dimension and the expected number of points in a fixed ball are of the same order of magnitude, and in fact asymptotically the same if $\rho_t / \ln (t) \to \infty$. 

We also mention that for a wide range of functionals, a shift in behaviour occurs when the limit of $\rho_t$ changes from zero to a nonzero constant to infinity. In our case, however, the shifts happen for different asymptotics of $\rho_t$ as seen in the theorem.

In the setting where the vertices arise from a binomial point process $X_n$ with parameter $n$ (as opposed to a Poisson point process), the clique number of $\c G(X_n,r_n)$ was investigated in \cite{M08},\cite{MM11}, the scan statistic in $\R^2$ in \cite{M99}, as well as a more general set of results in \cite{P02}.
In particular, \cite{M08},\cite{M99},\cite{P02} give two-point concentration results in the sparser regimes, while in \cite{MM11} an almost sure convergence is shown in all regimes, with the same asymptotics as in our results. 
Similarly to the clique number, the maximal degree also exhibits concentration on two consecutive values, discussed in \cite[Section~6]{PenroseBook}.
More recently, Kahle et al. \cite{KTW24} considered the clique number in perturbed (binomial) random geometric graphs.
We note that although the Poisson point process is conditionally binomial, we believe that the existing results for the binomial process are not granular enough to easily yield the corresponding statements for the Poisson case.
Additionally, we would like to mention that the results of \cite{ST12} that consider the maximum of certain score functions could theoretically be used to consider the dimension: $D_t^{\c V}\leq k$ is equivalent to the minimal diameter among $(k+1)$-tuples being at most $r_t$ ($D_t^{\c C}\leq k$ to the circumradius being at most $r_t/2$). However, though their results hold for arbitrary $k$, having $k=k_t\to \infty$ (for cases ii)-iv)) makes the $k$-dependent expressions that emerge in their computations complicated, and seemingly not tractable for our questions.

Our next result concerns large deviations, whose definition we shortly recall here (for more details see \cite{LDPbook}).
A family $(\PP_t)_{t>0}$ of probability measures on a topological space $\c X$ with $\sigma$-field $\c B$ is said to satisfy a large deviations principle (LDP) with speed $s_t\to \infty$ and good rate function $I:\c X\to [0,\infty]$ if the levels sets $\{x:I(x)\leq \alpha\}$ are compact for all $\alpha\in[0,\infty)$, and for all $B\in \c B$,
\begin{equation}\label{eq:LDP}
    -\inf_{x\in \inter B} I(x) \leq \liminf_{t\to \infty} s_t^{-1} \ln (\PP_t(B)) \leq \limsup_{t\to\infty} s_t^{-1} \ln (\PP_t(B))\leq - \inf_{x\in \cl B} I(x)
\end{equation}
holds, where $\inter(\cdot)$ and $\cl(\cdot)$ stand for the interior and closure operator, respectively. Further, we say that a family $(X_t)_{t>0}$ of random variables satisfies an LDP if their distributions do.

\begin{theorem}\label{thm:LDP}
Let $\eta_t$, $\rho_t$ and $D_t$ be as in \Cref{thm:main}, and assume $\rho_t/t^{-\alpha}\to \infty$ for all $\alpha>0$. Then $D_t/k_t$ satisfies an LDP with speed $m_t$ and good rate function $I(x)$, where
\begin{align*}
\intertext{i) for $\rho_t/\ln (t) \to 0$:}
k_t & = \displaystyle \frac{\ln (t) }{\ln (\ln (t) / \rho_t)},&  
m_t &= \ln (t) , & \text{ and } && 
I(x) &= \begin{cases} x-1, & \text{ if } x\geq 1, \\ \infty, & \text{ otherwise;}\end{cases} \\
\intertext{ii) for $\rho_t/\ln (t) \to B\in (0,\infty)$:}
k_t & = \beta  \cdot \frac{\kappa_d}{2^d} \rho_t & 
m_t & = \frac{\kappa_d}{2^d} \rho_t, & \text{ and } && 
I(x) &= \begin{cases} H(\beta x)-H(\beta), & \text{ if } x\geq 1, \\ \infty, & \text{ otherwise;}\end{cases}\\
\intertext{iii) for $\rho_t / \ln (t) \to \infty$:}
k_t & = \frac{\kappa_d}{2^d} \rho_t, & 
m_t &= \frac{\kappa_d}{2^d} \rho_t,& \text{ and } && 
I(x) &= \begin{cases} H(x), & \text{ if } x\geq 1, \\ \infty, & \text{ otherwise,}\end{cases}
\end{align*}
where $\kappa_d$ denotes the volume of the $d$-dimensional unit ball, $H(x)=1+x\ln x - x$, and $\beta>1$ is the unique number satisfying $H(\beta)=(B\cdot \kappa_d 2^{-d})^{-1}$.
\end{theorem}
To the author's knowledge, there are no precise distribution results known about the dimension in the $\rho_t /\ln (t) \to \infty$ regime. This is the topic of the following result for $d=1$.

\begin{theorem}\label{thm:Gumbel}
    Let $\eta_t$ be a homogeneous Poisson point process on $W=[0,1]\subset \R$, and let $\rho_t / \ln (t) \to \infty$. Then $D_t$ satisfies a Gumbel limit theorem; more precisely, 
    \[\PP\l(\frac{D_t-a_t}{b_t} \leq x\r) \to e^{-e^{-x}} \text{ as } t \to \infty,\] 
    where 
    \begin{equation}\label{eq:standardise}
    \begin{gathered}
        a_t =  \rho_t \l( 1+ \sqrt{\frac{2\ln(t/\rho_t)}{\rho_t}} \cdot \l(1+\frac{\ln(\ln(t/\rho_t))}{4\ln (t/\rho_t)} - \frac{\ln \sqrt{\pi}}{2\ln(t/\rho_t)}\r)\r) \\
        \text{ and } b_t = \sqrt{\frac{\rho_t}{2\ln(t/\rho_t)}}.
    \end{gathered}
    \end{equation}
\end{theorem}

This is obtained as a consequence of the following.

\begin{theorem}\label{thm:Poisson}
    Let $\eta_t$ be a homogeneous Poisson point process on $W=[0,1]\subset \R$, and let $\rho_t / \ln (t) \to \infty$.
    Let $I_n\coloneq [nr_t,(n+1)r_t)\cap [0,1]$ for $n=0,\ldots, N \coloneq \lfloor 1/r_t \rfloor$. For arbitrary $k\in \N$ and $n\in \{1,\ldots, N\}$, denote by $X_n^{(k)}$ the indicator of the event that $I_{n-1}\cup I_{n}$ contains an interval of length $r_t$ with at least $k$ points. Define $X^{(k)}\coloneq \sum_{n=1}^N X_n^{(k)}$, and let $Z^{(k)}\sim \Po(\EE X^{(k)})$. Choose $k=k_t(x)\in \N$ such that $k_t = a_t + xb_t + o(b_t)$ as $t\to \infty$, with $a_t,b_t$ given in \eqref{eq:standardise} and $x\in \R$ arbitrary. Then $\EE X^{(k_t(x))}\to e^{-x}$ as $t\to \infty$, and
    \[\dTV\big(X^{(k_t(x))},Z^{(k_t(x))}\big)= O\l(\frac{1}{\ln(1/r_t)}\r),\]
    where $\dTV$ denotes the total variation distance.
    In particular, convergence in distribution to a Poisson random variable with parameter $e^{-x}$ follows.
\end{theorem}

Clearly, we have $D_t\geq k-1 \Leftrightarrow X^{(k)}\geq 1$, and \Cref{thm:Gumbel} automatically follows.

Note that \Cref{thm:Gumbel} quantifies how much larger the dimension is compared to the number of points in an arbitrary fixed interval of length $r_t$, which is a lower bound on $D_t$. In particular, the latter is simply a Poisson random variable with parameter $\rho_t$, thus has fluctuations of order $\sqrt \rho_t$ around $\rho_t$. Compared to $\rho_t + O(\sqrt \rho_t$), the second-order term of the centralising function $a_t$ of $D_t$ has an additional logarithmic factor.

We suspect that a similar result to that of \Cref{thm:Gumbel} should hold in $\R^d$ for arbitrary $d>1$ as well, with the intervals being replaced by cubes, and as a consequence, a limit theorem for the dimension in the $l_{\infty}$-norm; however, getting such a precise result on the distribution appears to be difficult. 
This is because in $d=1$, the proof includes a purely combinatorial argument that relies on $\R^1$ having a natural linear order associated with it. Higher-dimensional spaces intrinsically do not possess such an order, and it appears that the idea cannot be easily extended.

Apart from the missing higher-dimensional setting, one could consider changing the underlying process as a potential extension. For a nonhomogeneous Poisson point process, we believe results corresponding to Theorems~\ref{thm:main} and \ref{thm:LDP} could be shown, following the tools used here and in \cite{MM11}. For \Cref{thm:Poisson} however, though we believe a similar result could hold, the nonhomogeneous setting would destroy the advantage gained from the combinatorial argument mentioned above, which relies heavily on uniformity. 
It would also be interesting to investigate the problem in substantially different settings as well, such as for Gibbs processes, though we expect the necessary methods to differ significantly from ours.

The paper is structured as follows. In \Cref{section:prelim}, we define some of the necessary objects in more detail, and recall some well-known results and computational tools. The proof of Theorems~\ref{thm:main} and \ref{thm:LDP} are split between Sections~\ref{section:exp} and \ref{section:tpc}: the former contains bounds on the tails in the regimes with $\rho_t$ asymptotically larger than all $t^{-\alpha}$ ($\alpha>0$), which are then used to prove the convergence in probability, expectation, and LDP; the latter section discusses the two-point concentration result, and the expectation in the regime where $\rho_t = O(t^{-\alpha})$ for some $\alpha>0$. Lastly, \Cref{section:dense} contains the proof of \Cref{thm:Poisson}.

\section{Preliminaries}\label{section:prelim}

%{\color{red} Poisson point process intro...} $\c N A$ denotes the number of points in $A$. Furthermore, let $\rho_t\coloneq tr_t^d$. 

We work in Euclidean space $\R^d$ with $d\geq 1$. The Euclidean norm is denoted by $\Vert\cdot \Vert$, and we write $\Vol(\cdot)$ for the volume (i.e. Lebesgue measure) of a measurable set.
We write $B^d(x,r)$ for the closed ball with center $x$ and radius $r$, and as a special case, $B^d \coloneq B^d(0,1)$ for the origin-centered unit ball, and $\kappa_d\coloneq \Vol (B^d)$ for its volume. We say that $K\subset \R^d$ is a convex body if it is convex, compact, and has nonempty interior.
%As the dimension $d$ will be fixed, we will mainly omit the superscript $d$ to simplify the notation.

We use the following asymptotic notation: % (Bachmann-Landau symbols): 
for arbitrary function $g$ and nonnegative $h$, we write
\begin{align*}
    g  = o(h)  \quad \text{ for } \quad &\lim_{x\to \infty} \frac{g(x)}{h(x)}  = 0; &
    g =\omega(h)  \quad \text{ for }\quad& \lim_{x\to \infty}\frac{g(x)}{h(x)}  = \infty;\\
    g  =O(h)  \quad \text{ for }\quad&\limsup_{x\to \infty} \frac{|g(x)|}{h(x)}<\infty; &
    g  =\Omega(h)  \quad \text{ for }\quad&\liminf_{x\to \infty} \frac{|g(x)|}{h(x)}>0;
\end{align*}
as well as $g=\Theta(h)$ when both $g=O(h)$ and $h=O(g)$, and $g\sim h$ for $\lim\limits_{x\to \infty} g(x) /h(x)=1$.
To avoid clutter throughout the proofs, we also write $\lesssim$, $\gtrsim$, and $\approx$ to denote $O(\cdot)$, $\Omega(\cdot)$, and $\Theta(\cdot)$, respectively.

To measure the closeness of distributions of real-valued random variables, we use the total variation distance $\dTV$ given by
\[\dTV(X,Y) \coloneq \sup_{A\in\c B(\R)} |\PP(X\in A)-\PP(Y\in A)|.\]

%We use the following asymptotic notation: % (Bachmann-Landau symbols): 
%for nonnegative functions $g$ and $h$, we say that $g=o(h)$ if $\lim_{x \to \infty} g(x)/h(x) = 0$; $g=O(h)$ if $ \limsup_{x\to \infty} g(x)/h(x) <\infty$; $g=\omega(h)$ if $\lim_{x\to \infty} g(x)/h(x)=\infty$, and $g = \Omega(h)$ if $\liminf_{x\to \infty} g(x)/h(x)>0$. We also write $g\ll h$ for $g=o(h)$. Furthermore, if $g=O(h)$ and $h=O(g)$, we say that $g=\Theta(h)$, and if $g(x)/h(x)\to 1$ as $x\to \infty$, we write $g\sim h$.

Let $W\subset \R^d$ be a fixed convex body, and assume $\Vol (W)=1$. Let $\mathbf N(W)$ denote the space of simple finite counting measures on $W$; note that generally, we will identify a counting measure with its support, that is, consider it as a finite set of points lying in $W$. Then the \textit{homogeneous Poisson point process} $\eta_t$ on $W$ with parameter $t>0$ is a random variable with values in $\mathbf N(W)$, satisfying the property that for a Borel set $A\subset W$, the number of random points in $A$, denoted by $\eta_t(A)$, is Poisson distributed with parameter $t\cdot \Vol(A)$, and the number of points in disjoint sets $A_1,\ldots, A_k\subset W$ are independent. We call $t$ the \textit{intensity} of the process. For a detailed introduction to the Poisson point process, we refer the reader to the book \cite{PoissonBook}.

%We also use the terminology of $U$-functionals: we say that $F = F(\eta_t)$ is a \textit{Poisson $U$-functional} with order $k\geq 1$ and kernel $f:(\R^d)^k\to \R$ if \[F= \sum_{(x_1,\ldots, x_k)\in \eta_{t,\neq}^k} f(x_1,\ldots, x_k),\]
%where $\eta_{t,\neq}^k$ denotes the set of distinct $k$-tuples of $\eta_t$.

An important property of the Poisson point process is the Mecke formula (see \cite[Theorem~4.4]{PoissonBook}). In its multivariate form, it states the following: if $f:(\R^d)^k\times \mathbf N(W) \to \R$ is a nonnegative function, then
\begin{equation}\label{eq:Mecke}
    \EE \sum_{(x_1,\ldots, x_k) \in \eta_{t,\neq}^k} f(x_1,\ldots, x_k,\eta_t) = t^k \int\limits_{W^k} \EE f(x_1,\ldots, x_k, \eta_t+\delta_{x_1} + \ldots + \delta_{x_k}) \dd x_1 \ldots \dd x_k,
\end{equation}
% \begin{multline*}
%     \EE \sum_{(x_1,\ldots, x_k) \in \eta_{t,\neq}^k} f(x_1,\ldots, x_k,\eta_t) =  \\ = t^k \int\limits_{W^k} \EE f(x_1,\ldots, x_k, \eta_t+\delta_{x_1} + \ldots + \delta_{x_k}) \dd x_1 \ldots \dd x_k,
% \end{multline*}
where $\eta_{t,\neq}^k$ denotes the set of distinct $k$-tuples of $\eta_t$.
%\[\EE \sum_{\substack{(x_1,\ldots, x_k) \\ \in \eta_{t,\neq}^k}} f(x_1,\ldots, x_k,\eta_t) = t^k \int\limits_{W^k} \EE f(x_1,\ldots, x_k, \eta_t+\delta_{x_1} + \ldots + \delta_{x_k}) \dd x_1 \ldots \dd x_k.\]
Note that if $f$ does not depend on $\eta_t$, then the expectation on the right hand side can be dropped.

In what follows, we consider the Vietoris--Rips and \v Cech complexes with vertex set $\eta_t$.
%We introduce the following indicator functions: let $\Delta^{n,\c V}_\delta = \Delta^{n,\c V}_\delta(x_1,\ldots, x_n)=1$ if $\Vert x_i-x_j\Vert \leq \delta$ for all $i,j=0,\ldots,k$, and $\Delta^{n,\c C}_\delta = \Delta^{n,\c C}_\delta(x_1,\ldots, x_n)=1$ if there exists a point $y\in \R^d$ that satisfies $\Vert x_i-y\Vert \leq \delta/2$ for all $i=0,\ldots, k$; set them to $0$ otherwise.
Recall that the tuple $\{x_0,\ldots, x_k\} \in \eta_{t,\neq}^{k+1}$ is in $\c V(\eta_t,r_t)$ if and only if their pairwise distances are at most $r_t$, and in $\c C(\eta_t,r_t)$ if and only if they are contained in a ball of radius $r_t/2$. The number of such $(k+1)$-tuples i.e. $k$-faces of the complex is denoted by $f_k^{\c V}$ and $f_k^{\c C}$ for the Vietoris--Rips and \v Cech complexes, respectively, where the omission of the superscript means either of the two constructions can be considered. In addition,  these tuples will sometimes be referred to as $k$-faces \textit{induced} by $\eta_t$; unless specified, this can mean both complexes.

Subgraph counts in the random geometric graph, and as a consequence, $f_k^{\c V}$ itself, is well-understood, and a large array of results is given e.g. in \cite[Section~3]{PenroseBook}. The computations therein can be easily adapted to $f_k^{\c C}$ of the \v Cech setting.
In particular, by \cite[Propositions~3.1 and 3.7]{PenroseBook},
\begin{equation}\label{eq:f-basics}
    \EE f_k \sim \frac{\mu_k}{(k+1)!} \cdot t \cdot \rho_t^k \quad \text{ and } \quad \VV f_k  \sim \sum_{i=1}^{k+1} \Phi_{k,i} \cdot t \cdot \rho_t^{2k-i+1} \text{ as } t \to \infty,
\end{equation}
holds for arbitrary $k\geq 0$,
where $\rho_t= tr_t^d$, and 
\[\mu_k\coloneq \int_{B^d} \ldots \int_{B^d} \indi\l(\{x_1,\ldots, x_k\} \text{ forms $(k-1)$-face with distance $1$}\r) \dd x_1 \ldots \dd x_k\]
for $k\geq 1$ and $\mu_0\coloneq 1$, and $\Phi_{k,i}$ is a positive constant depending only on $k$ and $i$.
Note that  these expressions differ for $k\geq 3$ for the two complexes, and so the constants in the asymptotic formulas do as well.
It is also known (see \cite[Theorem~3.4]{PenroseBook}) that if $\EE f_k$ tends to a nonzero constant $\lambda$, then $f_k$ tends in distribution to a Poisson distributed random variable with parameter $\lambda$.

We also state some standard bounds and asymptotics.
To bound binomial coefficients, we use the standard inequality $\binom nk\geq \l(\frac nk\r)^k$ for $n\geq k \geq 0$. 
For a Poisson-distributed random variable with parameter $\lambda$, denoted by $\Po(\lambda)$, the asymptotic
\begin{equation}\label{eq:PoissonAsy}
\PP(\Po(\lambda)=k)= e^{-\lambda} \cdot \frac{\lambda^k}{k!} \sim \frac{e^{-\lambda}}{\sqrt{2\pi k}}\cdot \l(\frac{\lambda e}{k}\r)^k
\end{equation}
holds as $k\to \infty$ by Stirling's formula. 
Observe that as a consequence, there exist universal constants $C>c>0$ such that 
\begin{equation}\label{eq:PoissonBounds}
c \cdot \frac{e^{-\lambda}}{\sqrt{k}}\cdot \l(\frac{\lambda e}{k}\r)^k \leq \PP(\Po(\lambda)=k)\leq C \cdot \frac{e^{-\lambda}}{\sqrt{k}}\cdot \l(\frac{\lambda e}{k}\r)^k
\end{equation}
for all $\lambda>0$ and $k\in \N$.
In addition, for the tails of a Poisson random variable, the following estimates, known as the Chernoff-Hoeffding bounds (see e.g. \cite[Lemma~1.2]{PenroseBook}), hold:
\begin{align}
    \PP(\Po(\lambda)\geq k) \leq e^{-\lambda}\cdot \l(\frac{e\lambda}{k}\r)^k & \text{ if } k\geq\lambda, \label{eq:PoissonUpper} \\
    \PP(\Po(\lambda)\leq k) \leq e^{-\lambda}\cdot \l(\frac{e\lambda}{k}\r)^k & \text{ if } k\leq\lambda. \label{eq:PoissonLower}
\end{align}

\section{Proof of the expectation and LDP}\label{section:exp}

In this section, we obtain some upper and lower bounds on the probability $\PP(D_t\geq k)$ for appropriately chosen $k$, in the regime where $\rho_t/t^{-\alpha} \to \infty$ as $t\to \infty$ for all $\alpha>0$.
First, we collect some technical details and preliminary results.

The first lemma is a standard statement in convex geometry. We present a short proof for completeness.
\begin{lemma} \label{lemma:balls}
    There are $\Theta(r_t^{-d})$-many disjoint translates of $r_t B^d$ contained in $W$.
\end{lemma}
\begin{proof}
    The upper bound is apparent from a volume comparison argument. 
    For a lower bound, consider the family of $r_t$-balls whose centers are given by the lattice $3r_t \Z^d$; they are disjoint by construction. 
    Any ball of fixed radius contains $\Omega(r_t^{-d})$-many of these balls, and since $W$ has nonempty interior, the statement follows.
    %Note that the lower bound of the lemma is a special case of \cite[Lemma~3.1]{M08}.
\end{proof}

The next lemma is the Poisson analogue of a statement that appears as part of a proof in \cite[Section~3.3,~Proof of the upper bound]{MM11}. We give a slightly simplified proof better suited to the Poisson setting, which was suggested by an anonymous referee.

Denote by $M_{A,r_t}\coloneq \max_{x\in \R^d} \eta_t(x+r_t A)$ the scan statistic with scanning set $r_t A$. 
%Since $\eta_t$ only contains finitely many points, this is indeed a maximum.

\begin{lemma}\label{lemma:scan}
    Let $A\subset \R^d$ be a convex body and $\eps>0$ arbitrary. Then for arbitrary integer $k\geq 0$,
    \[\PP(M_{A,r_t} \geq k) \leq  \frac{\Vol(r_t^{-1} W + (-A_{\eps}))}{1-e^{-\kappa_d \eps^d}} \cdot \PP(\Po(tr_t^d \Vol (A_\eps) )\geq k),\]
    where $A_\eps\coloneq A+\eps B^d$ denotes the $\eps$ parallel body of $A$.
\end{lemma}
Note that by the Steiner formula (see e.g. \cite[Theorem~6.6]{Gruber}),
\begin{equation}\label{eq:Steinerbound}
        \Vol(r_t^{-1} W + (-A_{\eps})) = r_t^{-d} \Vol(W) (1+O(r_t)) = r_t^{-d}(1+O(r_t))
\end{equation}
as $r_t\to 0$, where the constant implied by the big $O$ notation depends on $W$ and $A_{\eps}$.
\begin{proof}
%Since the scan statistic is defined as the maximum over a continuum of locations, an upper bound on its upper tail is obtained by relating it to an associated discrete problem.
    Let $\xi$ be a unit intensity Poisson point process on $\R^d$, independently from $\eta_t$. 

    For fixed points $p\in \R^d$ and $q\in B(p,r_t\eps)$, we have $p-q+r_tA \subset r_t\eps B^d + r_t A = r_t A_{\eps}$, or equivalently, $p+r_t A \subset q+r_t A_{\eps}$, thus $\eta_t(q+r_tA_\eps)\geq \eta_t(p+r_tA)$.
    %For fixed points $p\in \R^d$ and $q\in B(p, r_t\eps)$, we have $\eta_t(q+r_t A_\eps) \geq \eta_t(p+r_t A)$: this follows from the fact that $p-q+r_tA\subset r_t \eps B^d + r_t A = r_t A_{\eps}$, and so $p+r_t A \subset q+ r_t A_{\eps}$.    
    Further, for a fixed point $p\in \R^d$, the probability that $B(p,r_t\eps)$ contains a point of the rescaled point process $r_t\xi$ is exactly $1-e^{-\kappa_d \eps^d}$.
    These together imply that   
    \[\PP\l(\max_{x\in r_t\xi} \eta_t(x+r_t A_{\eps}) \geq k \ \bigg| \ \max_{x\in \R^d} \eta_t(x+r_t A) \geq k\r) \geq 1-e^{-\kappa_d \eps^d},\] 
    and consequently,
    \[\PP\l( M_{A,r_t}\geq k \r) \leq \frac{1}{1-e^{-\kappa_d \eps^d}} \PP\l(\max_{x\in \xi} \eta_t(r_t x + r_t A_\eps) \geq k \r).\] 
    
    The probability on the right-hand side can be upper bounded by
    \[ \PP\l(\max_{x\in \xi} \eta_t(r_tx + r_t A_\eps) \geq k \r) \leq \EE \l[\sum_{x\in \xi} \indi (\eta_t(r_t x + r_t A_\eps)\geq k)\r],\]
    and the bound evaluated through the Mecke formula \eqref{eq:Mecke}, yielding
    \begin{align*}
        \PP\l(\max_{x\in \xi} \eta_t(r_tx + r_t A_\eps) \geq k \r) 
            %& \leq \EE \l[\sum_{x\in \xi} \indi (\eta_t(r_t x + r_t A_\eps)\geq k)\r]\\
            & \leq \int_{\R^d} \PP(\eta_t(r_t x + r_t A_\eps)\geq k) \dd x \\
            & \leq \PP(\Po(tr_t^d \Vol(A_\eps))\geq k) \int_{\R^d} \indi ((r_t x + r_t A_\eps) \cap W \neq \emptyset) \dd x.
    \end{align*}
    The condition in the last integral is equivalent to $x\in r_t^{-1} W + (-A_\eps)$, and so the integral is precisely $\Vol(r_t^{-1} W + (-A_\eps))$.
\end{proof}
We now begin giving explicit tail bounds on $D_t$. Throughout the proof, let $\mu_t\coloneq \kappa_d 2^{-d} \rho_t = t\Vol((r_t/2) B^d)$.

\subsection{Proof of the lower bound.}

For the lower bound, the approach for the two constructions is the same. %, and so in this part $D_t$ denotes either $D_t^{\c C}$ or $D_t^{\c V}$. 
Also note that for higher moments, we obtain an automatic lower bound from Jensen's inequality, $\EE D_t^m \geq (\EE D_t)^m$, hence it is enough to consider the expectation.

For $\rho_t /\ln (t) \to \infty$, the lower bound of $\mu_t$ for $\EE D_t$ is trivial: for a fixed ball of radius $r_t/2$, its number of points is a lower bound for $D_t$, and thus $\EE D_t \geq \EE \Po(t\Vol((r_t/2) B^d)) = \mu_t$.
Otherwise, we use the following, equally straightforward geometric argument, also used in the proof of \cite[Lemma~4.1]{M08} and \cite[Lemma~3.8]{MM11}. If $D_t<n$, then every ball of radius $r_t/2$ contains at most $n$ points; further, by \Cref{lemma:balls}, there are $K =\Theta( r_t^{-d})$-many disjoint balls of radius $r_t/2$ contained in our observation window $W$. Hence it follows that \[\PP(D_t<n) \leq \PP(\Po(\mu_t)\leq n)^{K}.\]
For the $\mu_t /\ln (t) \to \infty$ case, the conclusion is simple, as $\PP(\Po(\mu_t)\leq \mu_t)\to 1/2$ follows from normal approximation, and thus $\PP(D_t < \mu_t)\to 0$. 
For the rest of the lower bound proof, assume $\mu_t = O(\ln (t))$.
Then we bound the expression above by
\[\PP(\Po(\mu_t)\leq n)\leq 1-\PP(\Po(\mu_t)=n+1) \leq \exp\{-\PP(\Po(\mu_t)=n+1)\},\] 
resulting in 
\[\PP(D_t<n)\leq \exp\{-K \cdot \PP(\Po(\mu_t) = n+1)\}= \exp\l\{-\Theta \l\{\PP(\Po(\mu_t) = n+1)\cdot \frac{t}{\mu_t} \r\}\r\},\]
where we also used that $K = \Theta(r_t^{-d}) = \Theta (t/\mu_t)$. Using \eqref{eq:PoissonBounds},
\begin{equation}\label{eq:alpha lower}
    - \ln \PP(D_t<n) \gtrsim \frac{t}{\mu_t} \cdot e^{-\mu_t} \cdot \l(\frac{\mu_t e}{n+1}\r)^{n+1} \cdot \frac{1}{\sqrt{n+1}} \eqcolon \alpha_n^{(L)}.
\end{equation}
%For the two remaining regimes,
Now, choose integers $(k_t)_{t>0}$ such that they satisfy the asymptotic behaviour stated in the theorem (for the appropriate regime). As all computations throughout this section are asymptotic as $t\to \infty$, we will see that the precise value of each $k_t$ is insignificant. 
We show that for arbitrary (but small enough) $\gamma>0$, we have with $n_t\sim k_t(1-\gamma)$ that $\alpha_{n_t}^{(L)}\to \infty$; the exact choice of $n_t$ is again irrelevant, as long as the asymptotic behaviour is met.
Showing $\alpha_{n_t}^{(L)}\to \infty$, implies that $\PP(D_t<k_t(1-\gamma))\to 0$, which implies the lower tail part of convergence in probability, as $\gamma$ was arbitrary. 
This is also sufficient for proving the lower bound on the expectation. By the Markov inequality, 
\[\frac{\EE D_t}{(1-\gamma)k_t} \geq \PP(D_t\geq (1-\gamma) k_t),\]
and thus proving that the right hand side probability converges to one yields that $\liminf_{t\to \infty} \EE D_t / k_t \geq 1-\gamma$, implying $\liminf_{t\to \infty} \EE D_t\geq k_t$ by the arbitrariness of $\gamma$.

For simplicity, consider
\begin{equation}\label{eq:beta lower}
  \beta_{n_t}^{(L)} \coloneq \ln (\alpha_{n_t}^{(L)}) = \ln (t)   -\ln (\mu_t)  - \frac 12 \ln (n_t) + (n_t-\mu_t) - n_t\ln \l(\frac{n_t}{\mu_t}\r)+ O(1).
\end{equation}
It is sufficient to show $\beta_{n_t}^{(L)}\to \infty$ as $t\to \infty$.
Observe that if $f(x)\to \infty$ as $x\to \infty$, then $f(x)\sim g(x)$ implies $\ln (f(x))\sim \ln (g(x))$ as $x\to \infty$. 
Consequently, any $(n_t)_{t>0}$ with the given asymptotic behaviour will yield the same $\beta_{n_t}^{(L)}$ asymptotically.

Assume first that $\mu_t \sim B' \ln (t)  = B\kappa_d2^{-d} \ln (t) $ for some $B>0$. Choose $\N \ni k_t\sim \beta \mu_t$, where $\beta>1$ is the unique number satisfying $\beta\ln \beta - \beta+1 = 1/B'$. 
Note that we then have $\ln (\mu_t) \sim \ln (\ln (t)) = o(\ln (t))$.
Plugging in $n_t \sim(1-\gamma)k_t = \beta(1-\gamma)\mu_t$, and the asymptotic behaviour of $\mu_t$, we have
\begin{align}\label{eq:beta mid}
    \beta_{n_t}^{(L)} 
    & = \ln (t)  - \frac 32 \ln (\mu_t) + (\beta(1-\gamma)-1)\mu_t - \beta(1-\gamma)\mu_t \ln (\beta(1-\gamma)) + O(1) \notag \\
    & = \ln (t)  \bigg[1 + o(1) + \beta(1-\gamma) B' - B' - \beta(1-\gamma)B' \ln (\beta (1-\gamma)) + o(1)\bigg].
\end{align}
After expanding and reorganising,
\[\beta_{n_t}^{(L)}  =  B' \ln (t) \l[\l(\frac{1}{B'}+ \beta - 1 - \beta \ln \beta \r) + \beta \bigg( \gamma \ln \beta - \gamma -  (1 - \gamma)\ln(1-\gamma)\bigg) + o(1)\r].\]
By choice of $\beta$, the expression in the first round parentheses is zero. For the second term, we have by the standard inequality $\ln(1-\gamma)\leq -\gamma$ that
\[\gamma \ln \beta - \gamma - (1-\gamma)\ln (1-\gamma) \geq \gamma \ln \beta - \gamma^2,\]
which is positive as long as $\gamma < \ln \beta$. Since $\beta>1$ is fixed and $\gamma$ is small, this may be assumed. This shows that $\beta_{n_t}^{(L)}$ is asymptotically some positive constant times $\ln (t) $, and hence tends to infinity. %, implying also $\alpha_{n_t}^{(L)}\to \infty$, which was our goal.

Assume now that $\mu_t /\ln (t) \to 0$, but $\mu_t/t^{-\alpha}\to \infty$ for all $\alpha>0$. Let $k_t \sim \ln (t)  / \ln (\ln (t) / \mu_t)$; by the choice of regime, the denominator is defined (for large enough $t$). In this regime, we additionally have $\ln (\mu_t)/\ln (t) \to 0$, as well as $k_t = o (\ln (t))$, since the denominator in its definition converges to infinity. Observe that the former expression can be negative ($\mu_t\to 0$, and so $\ln(\mu_t)\to -\infty$ is allowed), but the second condition on the regime ensures convergence to zero when divided by $\ln(t)$.
Plugging in $n_t \sim k_t(1-\gamma)$, we have
\begin{align*}
    \beta_{n_t}^{(L)}
     = \ln (t)  - &\ln (\mu_t)   - \frac 12 \ln \l(\frac{\ln (t) }{\ln (\ln (t) /\mu_t )}\r) +\frac{(1-\gamma)\ln (t) }{\ln (\ln (t) / \mu_t)} -\mu_t\\
     & -  \frac{(1-\gamma)\ln (t) }{\ln (\ln (t) / \mu_t)} \cdot \l[\ln \l(\frac{\ln (t) }{\mu_t}\r) + \ln \l( (1-\gamma) \ln \l(\frac{\ln (t) }{\mu_t}\r)\r)\r] + O(1),
\end{align*}
where the last expression in square parentheses is $\ln(n_t/\mu_t)$ broken apart at a convenient spot.
By the previous observations, the first line is exactly $\ln (t) (1+o(1))$, and the second line is $-(1-\gamma)\ln (t)(1 +o(1))$. In total, $\beta_{n_t}^{(L)} = \gamma\ln (t) (1+o(1))\to \infty$.
This concludes the proof of the lower bound.

\subsection{Proof of the upper bound.} 

%Let $\mu_t\coloneq \kappa_d 2^{-d} \rho_t$.

We have that the $m$-th moment can be written as \[\EE D_t^m = \sum_{i=1}^\infty i^m \PP(D_t = i) \leq  n^m + \sum_{i=n+1}^\infty i^m \PP(D_t \geq i)=: n^m + \alpha_n^{(U)}\]
for some integer $n$. We again show that for arbitrary $\gamma>0$, with $n=n_t\sim k_t(1+\gamma)$ we have that both $\PP(D_t \geq n_t)$ and $\alpha_{n_t}^{(U)}$ tend to $0$.

Proving the statement for the Vietoris--Rips and \v Cech complexes requires two somewhat different methods. 
In particular, we start with $D_t^{\c C}$, and later slightly modify the approach to be applicable  to $D_t^{\c V}$ as well.
This is because in the \v Cech complex, each simplex is constrained to an $r_t/2$-ball by definition. In contrast, a simplex of the Vietoris--Rips complex can a priori be more ‘spread out’: it is always contained in an $r_t$-ball, but generally not in an $r_t/2$-ball. Taking up more space in this sense could potentially be problematic when considering upper bounds, and thus this setting needs to be handled more delicately.

Considering $M_{A, r_t} = D_t^{\c C}+1$ with $A=\frac 12 B^d$, we can apply \Cref{lemma:scan} with some arbitrarily small fixed $\eps>0$. 
We then have $A_\eps = \frac 12 B^d + \eps B^d = \frac 12 (1+2\eps) B^d$ and $\Vol(A_\eps) = \kappa_d 2^{-d} (1+2\eps)^d$, yielding
\[\PP(D_t^{\c C}\geq i-1) \leq \frac{\Vol(r_t^{-1} W + (-A_\eps))}{1-e^{-\kappa_d \eps^d}} \cdot  \PP(\Po((1+2\eps)^d\mu_t)\geq i).\]
Using \eqref{eq:Steinerbound}, and the fact that $W$ and $\eps>0$ do not depend on $t$, the first factor is $O(r_t^{-d})=O(t/\mu_t)$.
Plugging this bound into the definition of $\alpha_n^{(U)}$, we have
\[\alpha_n^{(U)} = \sum_{i=n+2}^\infty (i-1)^m \PP(D_t\geq i-1)\lesssim \frac{t}{\mu_t} \sum_{i=n+2}^\infty i^m \PP(\Po((1+2\eps)^d \mu_t)\geq i).\]

By definition of a Poisson random variable, $N\PP(\Po(\lambda)=N) = \lambda \PP(\Po(\lambda)=N-1)$ for any $N\in \N$ and $\lambda>0$. It follows that for any $j\in \N$, %$\PP(\Po(\lambda)=N)\leq (\lambda/N)^j \PP(\Po(\lambda)=N-j)$
\begin{align*}
    \PP(\Po(\lambda)\geq N) 
    & = \sum_{i=N}^\infty\PP(\Po(\lambda)=i) \\
    & \leq \sum_{i=N}^\infty \l(\frac{\lambda}{N}\r)^j \PP(\Po(\lambda)=i-j) = \l(\frac{\lambda}{N}\r)^j \PP(\Po(\lambda)\geq N-j).
\end{align*}
As a consequence, 
\begin{align*}
    \alpha_{n}^{(U)}
    & \lesssim \frac{t}{\mu_t} \sum_{i=n}^\infty i^m \l(\frac{(1+2\eps)^d \mu_t}{i}\r)^{i-n+m}\PP(\Po((1+2\eps)^d \mu_t)\geq n-m) \\
    &\lesssim t\mu_t^{m-1}\PP(\Po((1+2\eps)^d \mu_t)\geq n-m)\sum_{i=n}^\infty \l(\frac{(1+2\eps)^d \mu_t}{n}\r)^{i-n}.
\end{align*}
If $n=n_t$ is chosen such that $\lim_{t\to \infty} (1+2\eps)^d \mu_t/ n_t <1$, the sum in the last expression is convergent, and does not depend on $t$ in the limit, i.e. can be upper bounded by a constant. 
Lastly, bounding the Poisson tail using \eqref{eq:PoissonUpper}, this yields that for such $n_t$,
\[\alpha_{n_t}^{(U)}\lesssim t\mu_t^{m-1} e^{-(1+2\eps)^d \mu_t} \l(\frac{e(1+2\eps)^d \mu_t}{n-m}\r)^{n-m}.\]
Now, fix $\gamma>0$, and choose $n_t\in \N$ such that the asymptotic $n_t\sim k_t(1+\gamma)$ is satisfied as $t\to \infty$; then $n_t-m \sim k_t(1+\gamma)$ holds as well, since $m$ is fixed. 
We choose $\eps=\eps(\gamma)>0$ so that the bound derived above hold; in particular we need $\lim_{t\to \infty} (1+2\eps)^d \mu_t/ ((1+\gamma)k_t) <1$. 
The details for each regime follow below.

Plugging in the asymptotic for $n_t$, we have
\begin{align}\label{eq:beta upper}
\begin{split}
    \beta_{n_t}^{(U)}\coloneq -\ln & (\alpha_{n_t}^{(U)}) 
    \gtrsim - \ln (t)  - (m-1)\ln (\mu_t) \\ 
    &+ (1+2\eps)^d \mu_t - (1+\gamma) k_t + (1+\gamma) k_t  \ln \l[\frac{k_t(1+\gamma)}{\mu_t (1+2\eps)^d}\r] + O(1).
\end{split}
\end{align} 
We show that this expression tends to infinity as $t\to \infty$ for arbitrarily small $\gamma$. This is turn yields $\alpha_{n_t}^{(U)}\to 0$, which is sufficient for proving the upper bounds of the theorem.

If $\mu_t/ \ln (t) \to \infty$, set $k_t=\mu_t$, and choose $\eps>0$ such that $(1+2\eps)^d = (1+\gamma)/2$. Then
\[\beta_{n_t}^{(U)} \gtrsim -\ln (t)  - (m-1) \ln (\mu_t) + \frac 12\mu_t (1+\gamma)\l[ 1-2 + 2\ln 2\r] + O(1).\]
The coefficient of $\mu_t$ is positive ($2\ln 2-1>0$ since $4>e$). By choice of regime, the $\mu_t$ term dominates the entire expression, which then converges to infinity.

If $\mu_t \sim B'\ln (t) =B\kappa_d 2^{-d} \ln (t) $, set $k_t=\beta \mu_t$ with $\beta>1$ satisfying $B'(\beta\ln \beta - \beta+1) = 1$, and $(1+2\eps)^d = 1+\gamma$, to get 
    \begin{align*}
        \beta_{n_t}^{(U)}  
        & \gtrsim - \ln (t)  - (m-1) \ln (\mu_t) + (1+\gamma)(1-\beta) \mu_t + (1+\gamma) \beta \mu_t \ln \beta + O(1)\\
        & = \ln (t)  \bigg[ - 1  - o(1)  + (1+\gamma) B'(1-\beta +\beta \ln \beta) + o(1) \bigg] \\
        & = \ln (t) (\gamma - o(1))  \to \infty, 
    \end{align*}
    where in the last step we used how we chose $\beta$. 

If $\mu_t / \ln (t) \to 0$, but $\mu_t/t^{\alpha}\to \infty$ for all $\alpha>0$, set $k_t=\ln (t)  /\ln(\ln (t) / \mu_t)$ and $(1+2\eps)^d = 1+\gamma$ to get 
    %\[\beta_{n_t}^{(U)} = \ln (t)  \l(-1+(1+\gamma)\l[1+\Theta\l(\frac{\ln \ln \frac{\ln (t) }{\mu_t}}{\ln \frac{\ln (t) }{\mu_t}}\r)\r] + \Theta\l(\frac{\ln \ln (t) }{\ln (t) }\r)\r)\to \infty.\]
\begin{align*}
    \beta_{n_t}^{(U)}
     \gtrsim - \ln (t)  - (m-1) & \ln (\mu_t)   + (1+\gamma)\mu_t - \frac{(1+\gamma)\ln (t) }{\ln (\ln (t) / \mu_t)} \\ 
     & + \frac{(1+\gamma)\ln (t) }{\ln (\ln (t) / \mu_t)} \cdot \l[\ln \l(\frac{\ln (t) }{\mu_t}\r) + \ln \l( \ln \l(\frac{\ln (t) }{\mu_t}\r)\r)\r].
\end{align*}
Summands two, three and four in the first lineare $o(\ln (t))$, so is the term coming from the second summand in the last expression in parentheses. We thus have
\[\beta_{n_t}^{(U)} \gtrsim \ln (t)  \l[ -1+(1+\gamma) - o(1)\r] = \ln (t)(\gamma - o(1))\to \infty.\]

This completes the proof for the \v Cech complex. For the Vietoris--Rips complex, we use the approach appearing in \cite{MM11}.

Fix $\eps'>0$, and let $A_1,\ldots, A_K\subset \eps' \Z^d$ be all the subsets of $\eps' \Z^d$ that contain the origin and have diameter at most $1 + 2c_0\eps'$, where $c_0\coloneq \diam ([0,1]^d)$. Also let $W_j\coloneq \conv{A_j}$ $(j=1,\ldots, K)$. 
Assume now that the points $X_0,\ldots, X_i$ form an $i$-face of $\mathcal V(\eta_t,r_t)$, meaning that each pairwise distance is at most $r_t$. Then the point set $Y_j\coloneq (X_j-X_0)\cdot r_t^{-1}$ for $j=0,1,\ldots i$ contains the origin and has diameter at most $1$.
Consider the unit lattice cubes of $\eps'\Z_d$ that contain at least one of the $Y_j$, and take the set $A\subset \eps'\Z^d$ of all of their vertices. These are exactly the lattice points having distance at most $c_0\eps'$ from some $Y_j$, and thus
\[A = \{p\in \eps' \Z^d \colon \Vert p-Y_j\Vert \leq c_0\eps'\text{ for some } j=0,\ldots, i\}.\]
The set $A$ then contains the origin (as $Y_0=0$) and has diameter at most $1+2c_0\eps'$, hence it coincides with one of the sets $A_j$ defined above, say $A_{j_0}$. By construction, each lattice cell that contains some $Y_j$ is in $\conv A$. Thus $Y_0,\ldots, Y_i\in \conv A$, and consequently, $X_0,\ldots,X_i \in X_0 + r_t \conv A$. This implies that $M_{W_{j_0},r_t}\geq i+1$, hence $\max_{j=1,\ldots,K} M_{W_j,r_t} \geq D_t^{\mathcal V}-1$.

Using this statement, as well as \Cref{lemma:scan}, we obtain that for arbitrary $\eps>0$, 
\begin{align*}
    \PP(D_t^{\mathcal V} \geq i-1) 
    &\leq \sum_{j=1}^K \PP(M_{W_j,r_t}\geq i) \\
    &\leq \sum_{j=1}^K \frac{\Vol(r_t^{-1} W + (-W_{j,\eps}))}{1-e^{-\kappa_d \eps^d}} \PP(\Po( \rho_t \Vol (W_{j,\eps}))\geq i),
\end{align*}
where $W_{j,\eps}\coloneq (W_j)_{\eps} = W_j + \eps B^d$. 
Observe that by construction, $\diam (W_{j,\eps})\leq 1+2c_0\eps'+ 2\eps$.
On the one hand, this implies that the constants in the first factor can be made universal. In particular, $\Vol(r_t^{-1}W + (-W_{j,\eps}))\leq \Vol(r_t^{-1}W+(1+c_0\eps'+2\eps)B^d)$, which is $r_t^{-d}(1+O(r_t))$ by \eqref{eq:Steinerbound}, with the implied constant not depending on $j$.
On the other hand, it is known from the isodiametric inequality (see e.g. \cite[Theorem~8.8]{Gruber}) that among convex bodies with fixed diameter, the ball has maximal volume, implying here that $\Vol (W_{j,\eps}) \leq \l(1+2c_0\eps'+2\eps\r)^d 2^{-d}\kappa_d$. 
Lastly, since $K=K(\eps')$ is a constant, we have
\[\PP(D_t^{\mathcal V} \geq i-1) \lesssim \frac{t}{\mu_t} \cdot \PP(\Po(\mu_t (1+2c_0\eps' +2\eps)^d)\geq i).\]
The proof can then be carried out the same way as for the \v Cech complex. \hfill $\qed$

\subsection{Proof of the LDP}

To prove \Cref{thm:LDP}, we simply use the bounds on the distribution obtained above. First, we show that the defining chain of inequalities \eqref{eq:LDP} holds for all half-intervals $[a,\infty)$, which then easily extends to all intervals $[a,b)$. Second, we show that the family of probability measures is exponentially tight, which, together with the statement for the intervals, yields the full LDP (see \cite[Theorem~4.1.11, Lemma~1.2.18]{LDPbook}).

Let $k_t$ be as stated in the theorem, according to the asymptotic behaviour of $\rho_t$. Recall also that $H(x)=1+x\ln x-x$.

For $a<1$, it was shown in the proof of \Cref{thm:main} that $\PP(D_t\geq ak_t)\to 1$, and so $\inf_{x>a} I(x)=0$ follows. (Alternatively, this follows since $D_t/k_t\to 1$ in probability, which we have already shown.)

For $a>1$, we have $-\ln \PP(D_t < ak_t)\gtrsim \alpha_{ak_t}^{(L)}$ as defined in \eqref{eq:alpha lower}. 
To find its asymptotic behaviour, plug $n_t = ak_t$ into \eqref{eq:beta lower}.
For $\mu_t/\ln (t)\to \infty$, this yields
\begin{align*}
    \beta_{ak_t}^{(L)} 
    = \ln(\alpha_{ak_t}^{(L)}) 
    & = \mu_t \l[\frac{\ln(t)}{\mu_t} - \frac{\ln(\mu_t)}{3\mu_t}+(a-1)-a\ln a + o(1) \r] \\
    & = -\mu_t \big[ H(a) + o(1)\big] \sim -m_t I(a).
\end{align*}
For $\mu_t\sim B'\ln(t) = B\kappa_d 2^{-d}\ln (t)$, plugging $1-\gamma=a$ into \eqref{eq:beta mid} gives 
\begin{align*}
    \beta_{ak_t}^{(L)} 
    & = \ln (t) \bigg [ 1 - B'\big( -a\beta + 1 + a\beta \ln(a\beta) \big)\bigg] \\ 
    &= -B'\ln(t) \big[H(a\beta)-H(\beta) \big] \sim -m_t I(a),
\end{align*}
using that $H(\beta)=1/B'$ by definition.
Lastly, for $\mu_t/\ln(t)\to 0$, we have shown before that $\beta_{ak_t}^{(L)} \sim (1-a)\ln (t) = -m_t I(a)$.
This additionally means that $\beta_{ak_t}^{(L)} \to -\infty$ as $t\to \infty$ in all three regimes, implying $\alpha_{ak_t}^{(L)}\to 0$ and $\PP(D_t < ak_t)\to 1$. 
Thus
\[\PP(D_t \geq  ak_t) = 1- \PP(D_t < ak_t) \sim -\ln \PP(D_t< a k_t) \gtrsim \alpha_{ak_t}^{(L)},\]
and hence 
\[1\geq
\limsup_{t\to \infty} \frac{\ln \PP(D_t \geq ak_t)}{\ln (\alpha_{ak_t}^{(L)})}  = 
\limsup_{t\to \infty} \frac{\ln\PP(D_t \geq ak_t)}{\beta_{ak_t}^{(L)}} =
\limsup_{t\to \infty} \frac{\ln \PP(D_t \geq  ak_t)}{-m_tI(a)}.\]
Observe that though there is a constant factor between $\PP(D_t\geq ak_t)$ and $\alpha_{ak_t}^{(L)}$ hidden in the $\gtrsim$ notation, taking the logarithm makes that term negligible.
Note also that setting $1-\gamma=a$ above requires $\gamma<0$, which is not what was considered for the expectation, but since \eqref{eq:beta lower} holds for arbitrary $n_t> 0$, any $\gamma$ with $1-\gamma> 0$ is allowed in the case distinctions that follow.

For the lower bound, we can show with a similar computation using \eqref{eq:beta upper} (set $1+\gamma=a$) that there is a constant $C=C(a)>0$ such that 
\[\lim_{t\to \infty} \frac{\beta_{ak_t}^{(U)}}{m_t} \geq I(a)(1-C\eps)\] 
for any small enough $\eps>0$. In particular, the bound holds whenever $\eps>0$ satisfies $(1+2\eps)^d<a$; since $a>1$, this is indeed true for small enough $\eps$.
Since $\PP(D_t\geq ak_t)\lesssim \alpha_{ak_t}^{(U)}$ by construction, this yields 
\[1\leq   
\liminf_{t\to\infty} \frac{\ln \PP(D_t \geq ak_t)}{\ln (\alpha_{ak_t}^{(U)})}
= \liminf_{t\to\infty} \frac{\ln \PP(D_t \geq ak_t)}{-\beta_{ak_t}^{(U)}}
\leq \liminf_{t\to \infty} \frac{\ln \PP(D_t \geq ak_t)}{-m_tI(a)(1-C\eps)}.\]
Since $\eps>0$ was arbitrary, this yields that $\liminf_{t\to \infty} \ln \PP(D_t\geq ak_t)/(-m_tI(a))\geq 1.$
Together with the limit superior, this shows that $\ln \PP(D_t\geq ak_t)\sim -m_tI(a)$ for all $a>1$.

From this, it easily follows that the defining inequalities hold for $[a,b)$ for arbitrary $a<b$. If $a<1$ or $b<1$, the conclusion is essentially trivial. For $a,b\geq 1$, we have
\[-\ln \PP(a\leq D_t /k_t < b) = -\ln \PP(D_t/k_t\geq a) - \ln \l(1-\frac{\PP(D_t/k_t\geq b)}{\PP(D_t/k_t\geq a)}\r).\]
The ratio of probabilities in the second logarithm tends to zero, since its logarithm $\ln \PP(D_t/k_t\geq b)-\ln \PP(D_t/k_t\geq a)\to -\infty$, following from the monotonicity of the rate function $I(x)$ (on the interval where $I(x)<\infty$).

Lastly, exponential tightness (see Equation (1.2.17) in \cite{LDPbook}) follows easily from $I(x)\to \infty$ $(x\to \infty)$: for arbitrary $M$, there is $K$ such that $I(K) > M$, and thus \[\limsup_{t\to \infty} \frac{\log \PP(D_t/k_t\in [0,K]^c)}{m_t} < -M,\]
with $\cdot^c$ denoting the complement. 
This concludes the proof of the full LDP. \hfill $\qed$

\section{Proof of the two-point concentration}\label{section:tpc}

In this section, we prove the two-point concentration results of \Cref{thm:main}, as well as the expectation in the first regime. 
%We give a quick and easy proof for the first regime for the sake of completeness (this is very similar to \cite{M99}*{Section~5.4}), and give the Poisson point process version of the proof carried out in \cite{M08} in the second regime.

First, let $\rho_t\coloneq tr_t^d = O(t^{-\alpha})$ for some $\alpha>0$. This case lends itself to a very straightforward proof using the known properties of the $f$-vector, or more precisely, the number of $n$-faces $f_n$; see also \cite[Theorem~6.1]{PenroseBook}. Namely, on one hand, we have for every integer $n\geq 0$ by \eqref{eq:f-basics} that $\EE f_n = \Theta (t\rho_t^n) = O(t^{1-n\alpha})$; on the other hand, we have $\rho_t\to 0$ as $t\to \infty$, and thus $\VV f_n = \Theta\l(t \rho_t^n\r)=O(t^{1-n\alpha})$ as well.
Then for $n>1/\alpha$, we have $\EE f_n\to 0$, and clearly $\EE f_0 = t \to \infty$. It follows that a unique integer $k$ exists such that $\EE f_{k-1} \to \infty$ and $\EE f_k\to \lambda<\infty$; clearly, $k\geq 1$. Note that by the monotonicity of $r_t$, $\EE f_k$ is indeed convergent, as well as that $\lambda$ might be zero or positive. Further, we have that $\EE f_n\to 0$ for all $n>k$.

%It is intuitively clear that $\EE f_{k-1}\to \infty$ implies that there is a $(k-1)$-face with high probability, and $\EE f_{k+1}\to 0$ implies that a $(k+1)$-face occurs with low probability,  hence $D_t$ is at least $k-1$ but less than $k+1$ with high probability. We make this argument precise as follows.

%\[\EE D_t^l = \sum_{k=1}^\infty \PP(D_t^l\geq k) = \sum_{n=1}^\infty (n^l-(n-1)^l)\PP(D_t\geq n) = \sum_{n=1}^\infty (n^l-(n-1)^l) \PP(f_n>0)\]

%Let $k_1 < k_2$ be arbitrary nonnegative integers. Then clearly $f_{k_2}\geq 1$ implies $f_{k_1}\geq \binom{k_2+1}{k_1+1}$, and by Markov's inequality and a well-known bound on binomial coefficients, we obtain \[\PP(f_{k_2}\geq 1) \leq \PP\l(f_{k_1} \geq \binom{k_2+1}{k_1+1} \r)\leq \frac{\EE f_{k_1}}{\binom{k_2+1}{k_1+1}}\leq \EE f_{k_1} \cdot \l(\frac{k_1+1}{k_2+1}\r)^{k_1+1}.\]

%Let $K$ be the greatest index $k$ for which $\EE f_k$ tend to infinity as $t\to\infty$, i.e. $\EE f_k \to \infty$ as $t\to\infty$ for all $k\leq K$, while $\EE f_k$ tends to a finite constant for every $k>K$, and let $K'$ be the smallest index $k$ with $\EE f_k \to 0$. Note that either $K'=K+1$ or $K'=K+2$ holds: the $\lim_{t \to \infty} \EE f_{K+1} =: \alpha$ can be either zero or positive, but the subsequent expectations tend to $0$ (...).

%First, note that if $\alpha>0$, it is known REF that $f_k$ tends in distribution to a Poisson random variable with parameter $\lambda$, hence $\PP(f_k>0) \to 1 - e^{-\lambda}$.

By Markov's inequality, 
\[\PP(f_n\geq 1) \leq \EE f_n\] 
and 
\[\PP(f_n = 0) = \PP(f_n\leq 0) =\PP(\EE f_n - f_n \geq \EE f_n) \leq \frac{\VV f_n}{(\EE f_n)^2} = \Theta\l( \frac{1}{\EE f_n}\r)\]
for arbitrary $n$.
As a consequence, we immediately obtain that 
\begin{align*}
    \PP(D_t>k) & = \PP(f_{k+1} \geq 1) = O\big(t\rho_t^{k+1}\big) \\ %\text{ and }
    \PP(D_t<k-1)& = \PP(f_{k-1}=0) = O\big((t\rho_t^{k-1})^{-1}\big),
\end{align*}
with both expressions tending to $0$ as $t\to \infty$.
%As $\rho_t\to 0$, \eqref{eq:f-basics}) implies that $\EE f_n \sim \VV f_n$ as $t\to \infty$, hence $\VV f_n / (\EE f_n)^2 \sim 1/\EE f_n$,
%\[\lim_{t\to \infty} \PP(D_t<k-1) = \lim_{t\to \infty} \PP(f_k=0) \leq \lim_{t\to \infty} \frac{\VV f_k}{(\EE f_k)^2} = \lim_{t\to\infty} \frac{1}{\EE f_k} = 0\]
Furthermore, if $\lambda >0$, then $f_k$ is asymptotically Poisson distributed with parameter $\lambda$. More precisely, by \cite[Theorem~3.4]{PenroseBook}, $\dTV(f_k, \Po(\EE f_k)) = O\big(t\rho_t^{k+1}\big)$, and thus 
\begin{align*}
    |\PP(D_t\geq k)-(1-e^{-\lambda})|
    & =|\PP(f_k\geq 1)-(1-e^{-\lambda})| \\ 
    & = O\big(t\rho_t^{k+1}\big)+ |e^{-\lambda}-e^{-\EE f_k}| = O\big(t\rho_t^{k+1}\big) + O(\EE f_k-\lambda).
\end{align*}
If $\lambda=0$, we have $\PP(f_k\geq 1) = O(t\rho_t^k) = O(\EE f_k-\lambda)$ formally.

Using that
\[\dTV(D_t,Z)= \frac 12 \sum_{n=1}^\infty |\PP(D_t=n)-\PP(Z=n)|,\]
the statement for the total variation distance follows.

Lastly for the $\rho_t=O(t^{-\alpha})$ regime, we consider the moments. Let $l\coloneq \max\{m,k\}$. It follows from the previous argument that 
\begin{multline*}
    \l| \sum_{n=0}^{l} n^m\PP(D_t=n) -((k-1)^m e^{-\lambda} + k^m (1-e^{-\lambda}))\r| \\ = O\big( (t\rho_t^{k-1})^{-1}\big)+O\big( t\rho_t^{k+1} \big) + O(\EE f_k-\lambda).
\end{multline*}  %Then it is again clear that $\sum_{n=k+1}^{l} n^m \PP(D_t=n) \to 0$, since it is a finite (possibly empty) sum where each term tends to $0$ as $t\to \infty$. 
Hence it remains to be shown that the tail starting at $n=l+1$ has at most this order.
By definition, any $(l+2)$-element subset of an $n$-face is an $(l+1)$-face, and so if there is at least one $n$-face, there are at least $\binom{n+1}{l+2}$-many $(l+1)$-faces. As a consequence, by Markov's inequality, as well as bounds on the binomial coefficient, \[\PP(f_n\geq 1)\leq \PP\l(f_{l+1}\geq \binom{n+1}{l+2}\r)\leq \frac{\EE f_{l+1}}{\binom{n+1}{l+2}} \leq \EE f_{l+1} \cdot \l(\frac{l+2}{n+1}\r)^{l+2},\]
for all $n\geq l+1$, and it follows that \[\sum_{n=l+1}^\infty n^m \PP(D_t\geq n) \leq  (l+2)^{l+2} \EE f_{l+1} \sum_{n=l+1}^\infty \frac{1}{(n+1)^{l+2-m}}.\]
The sum is convergent, since $l+2-m\geq 2$ as a consequence of $l\geq m$, and $\EE f_{l+1} = \Theta(t\rho_t^{l+1})=O(t\rho_t^{k+1})$ by $l\geq k$, which together imply the statement.

We now turn to the proof of the other regime, where $\rho_t/t^{-\alpha}\to \infty$ for all $\alpha>0$, but $\rho_t / \ln (t) \to 0$, where we again prove a two-point concentration, but with two consecutive integers that depend on $t$. Unfortunately, the number of $n$-faces is not a good (enough) indicator for $D_t$ anymore. The following argument comes from \cite{M08}, altered to fit the Poisson setting instead of the binomial therein.

Let $N_n$ denote the number of points that are part of some $n$-face, and $M_n$ the number of pairs of points that are contained in some $n$-face (together). On one hand, $D_t<n$ is equivalent to $N_n=0$; on the other hand, $D_t\geq n$ implies that $M_n\geq \binom{n+1}{2}$. Thus by Markov's inequality, we have
\begin{align}\label{eq:midtpcI}
    \begin{split}
        \PP(D_t<n) & =\PP(N_n=0) \leq \frac{\VV N_n}{(\EE N_n)^2}  \quad \text{ and } \\ \PP(D_t\geq n )&\leq \PP\l(M_n\geq \binom{n+1}{2}\r) \leq \frac{2\EE M_n}{n(n+1)}
    \end{split}
\end{align}

It seems that the reason these functionals describe the behaviour of the dimension better than $f_n$ is that for large $n$ (in particular, $n=n_t\to \infty$), even a single $n+k$-face supplies a large number of $n$-faces (in particular, $\binom{n+k+1}{n+1}$), which is large in $n$ even for fixed $k$; thus the expected number of $n$-faces can be large, even if their existence is not very likely. Considering instead the {\it number of points} in $n$-faces (somewhat) gets rid of this bias.

For a point configuration $\zeta\subset \R^d$ and $x\in \zeta$, write $f_n(x,\zeta)$ and $\c F_n (x,\zeta)$, respectively, for the number and set of $n$-faces induced by $\zeta$ that $x$ is contained in. 

We start with $N_n$ and the lower tail. Applying the Mecke formula, we have
\begin{align*}
    \EE N_n & = \EE \sum_{x\in \eta_t} \indi\l(f_n(x,\eta_t)\geq 1\r) = t \int_W \PP(f_n(x,\eta_t\cup\{x\})\geq 1)\dd x, \\
    \EE N_n^2 &= t\int_W \PP(f_n
    \begin{multlined}[t]
        (x,\eta_t\cup\{x\}) \geq 1)\dd x \\ 
        + t^2 \iint_{W^2} \PP(f_n(x,\eta_t\cup\{x,y\}),f_n(y,\eta_t\cup\{x,y\})\geq 1) \dd x \dd y,
    \end{multlined}
\end{align*}
where for the second moment we considered $N_n^2$ as a double sum over points $x$ and $y$, dividing said sum into $x=y$ and $x\neq y$ parts, to which the Mecke formula can be applied (separately).
%The expression is then divided into two regions: pairs of points with distance above or below $2r_t$.

The expression $\PP(f_n(x,\eta_t\cup\{x\})\geq 1)$ is the same for all $x$ for which $B(x,r_t)\subset W$; we denote the set of such points (called the inner parallel body of $W$) by $W_{-r_t}$, and write $p_{n,t}$ for this common probability. In addition, for points $z$ outside this set, we have $\PP(f_n(z,\eta_t\cup\{z\})\geq 1)\leq p_{n,t}$. As a consequence,
\[tp_{n,t} \geq \EE N_n \geq \Vol(W_{-r_t})\cdot tp_{n,t}.\]
Further, by a standard convex geometry result (see \cite[Section~6]{Gruber}), $\Vol(W_{-r_t})= \Vol (W) - O(r_t) = 1- O(r_t)$, and we have that $\EE N_n \sim tp_{n,t}$.

Now, we evaluate the second summand of the second moment expression above. We upper bound it by the sum of the following three integrals:
\begin{align*}
    I_1 & \coloneq t^2 \int_W \int_{W\setminus B(x,2r_t)} \PP(f_n(x,\eta_t\cup\{x,y\}),f_n(y,\eta_t\cup\{x,y\})\geq 1) \dd y \dd x, \\
    I_2 & \coloneq t^2 \int_W \int_{B(x,2r_t)\cap W} \PP(f_n(x,\eta_t\cup\{x\})\geq 1) \dd y \dd x, \\
    I_3 & \coloneq t^2 \int_W \int_{B(x,2r_t)\cap W} \PP(f_n(x,\eta_t\cup\{x,y\})\geq 1,\  f_n(x,\eta_t\cup\{x\})=0) \dd y \dd x.
\end{align*}
This indeed gives an upper bound: first we divide the $y$-domain into two parts based on distance to $x$, then upper bound the probability in the original integral by the sum of the probabilities in $I_2$ and $I_3$.

If $\Vert x-y\Vert > 2r_t$, then $B(x,r_t)\cap B(y,r_t) = \emptyset$, and hence the regions where their respective neighbours are located are disjoint. As a consequence, they are contained in an $n$-face independently from each other, and in particular we have
\[I_1\leq t^2 \int_W\PP(f_n(x,\eta_t\cup\{x\})\geq 1) \int_{W\setminus B(x,2r_t)}  \PP(f_n(y,\eta_t\cup\{y\})\geq 1)\dd y \dd x  \leq (\EE N_n)^2,\]
where in the second inequality we extended the second integral to $W$, thus getting two independent integrals each corresponding to the expectation.

To evaluate $I_2$, we can integrate with respect to $y$, and obtain
\[I_2\leq t^2 2^d \kappa_d r_t^d \int_W \PP(f_n(x,\eta_t\cup\{x\})\geq 1)\dd x = 2^d\kappa_d\rho_t \EE N_n.\]
Lastly, consider $I_3$.
The event in its definition implies that $x$ and $y$ must form an edge, i.e. $y\in B(x,r_t)$, and $y$ must be contained in all $n$-faces in $\c F_n(x, \eta_t \cup \{x,y\})$. 
Also note that by construction, $\c F_n(x,\zeta)=\c F_n(x,\zeta\cap B(x,r_t))$ for any point configuration $\zeta\subset \R^d$. 
Writing $\xi_{t,x} \coloneq \eta_t\cap B(x,r_t)$ for $x\in W$, this yields
\begin{equation*}
I_3\leq t^2\int_W \int_{B(x,r_t)\cap W} 
    \PP \bigg(f_{n}\big(x, \xi_{t,x}
    \begin{multlined}[t]
    \cup\{x,y\}\big )\geq 1,  \\
    y\in \bigcap \bigg\{ \sigma \in \c F_n  \big(x,\xi_{t,x}\cup\{x,y\}\big) \bigg\} \bigg)\dd y \dd x.
    \end{multlined}
\end{equation*}
Write $v_x \coloneq \Vol(B(x,r_t)\cap W)$, and let $Y, X_1, X_2,\ldots$ be independent uniform points in $B(x,r_t)\cap W$. Condition on the number of points in $\xi_{t,x}$ to obtain
\[I_3\leq t^2 \int_W \sum_{N=n-1}^\infty J_{N,x} \PP(\Po(tv_x)=N)\dd x\]
with
\begin{align*}
    J_{N,x}
    & \coloneq \int_{B(x,r_t)\cap W} \PP \bigg(f_{n} \big(x, \{x,X_1,
    \begin{multlined}[t]
    \ldots, X_N,y\}\big )\geq 1,  \\
    y\in \bigcap \bigg\{ \sigma \in \c F_n  \big(x,\{x,X_1,\ldots, X_N,y\}\big) \bigg\} \bigg)\dd y
    \end{multlined}\\
    & = v_x \PP \bigg(f_{n} \big( x, \{x, X_1,
    \ldots, X_N,
    \begin{multlined}[t]
    Y\}\big )\geq 1,  \\
    Y\in \bigcap \bigg\{ \sigma \in \c F_n  \big(x,\{x, X_1,\ldots, X_N,Y\}\big) \bigg\} \bigg).
    \end{multlined}
\end{align*}
The last expression is now symmetric in $X_1,\ldots, X_N,Y$. For each realisation, at most $n$ of these $N+1$ points are contained in all $n$-faces considered ($n$, and not $n+1$, since $x$ is contained in all of them by construction). 
As a consequence, we have
\[J_{N,x} \leq \frac{v_x\cdot n}{N+1} \cdot \PP\l(f_{n} (x,\{x,X_1,\ldots, X_{N+1}\}) \geq 1\r).\]
Plugging this into the bound on $I_3$, we have
\begin{align*}
    I_3
    & \leq tn \int_W \sum_{N=n-1}^\infty \PP(f_{n}(x,\{x,X_1,\ldots, X_{N+1}\})\geq 1) \cdot \frac{tv_x}{N+1} \PP(\Po(tv_x)=N) \dd x \\
    & =  tn \int_W \sum_{N=n-1}^\infty \PP(f_{n}(x,\{x,X_1,\ldots, X_{N+1}\})\geq 1) \cdot \PP(\Po(tv_x)=N+1) \dd x \\
    & \leq tn \int_W \PP(f_n(x,\xi_{t,x}\cup\{x\})\geq 1 ) \dd x = n \EE N_n.
\end{align*}
Putting together the three bounds, we have $\VV N_n \leq(1+2^d\kappa_d\rho_t+n)\EE N_n$. 

For the upper tail in \eqref{eq:midtpcI}, we have by the Mecke formula that
\[\EE M_n = \frac{1}{2} t^2 \int_{W^2} \PP(\c F_n(x,\eta_t\cup\{x,y\})\cap \c F_n(y,\eta_t\cup\{x,y\})\neq \emptyset)\dd x \dd y.\]
%\[\EE M_n = t^2 \int_{W^2} \PP(\text{$x$ and $y$ are in an $n$-face induced by $\eta_t+\delta_x+\delta_y$})\dd x \dd y.\]
For this event to hold, there must be an edge between $x$ and $y$, and $\c F_{n-1}(x,\eta_t\cup\{x\})$ must be nonempty. Hence
\begin{align*}
    \EE M_n
    & \leq \frac 12 t^2 \int_W \int_{B(x,r_t)} \PP(f_{n-1}(x,\eta_t\cup\{x\})\geq 1)\dd y \dd x \\
    & \leq \frac{\kappa_d}{2} t r_t^d \EE N_{n-1} = \frac{\kappa_d}{2}\rho_t \EE N_{n-1}.
\end{align*}
Plugging these bounds into \eqref{eq:midtpcI}, we have
\begin{equation}\label{eq:midtpcII}
    \PP(D_t<n-1) \leq \frac{n+2^d\kappa_d\rho_t}{\EE N_{n-1}}\ \text{ and }\ \PP(D_t>n) \leq \frac{\kappa_d \rho_t}{n}\cdot \frac{\EE N_n}{n+1}.
\end{equation}
We aim to find $n$ such that both expressions tend to zero as $t\to \infty$.

Let $m_t \coloneq \ln (t) /\ln (\ln (t) / \rho_t)$; note that $m_t/\rho_t\to \infty$. It is easy to see that the function $\EE N_n - (n+1) \cdot (m_t/\rho_t)^{1/2}$ is monotone decreasing in $n$ (for fixed $t$), positive for $n=0$ and negative for any $n=n_t$ with $n_t\sim t$. Thus there is a (unique) integer $k=k_t$ such that
\[\EE N_{k_t-1}-k_t\l(\frac{m_t}{\rho_t}\r)^{1/2} \geq 0 > \EE N_{k_t}-(k_t+1)\l(\frac{m_t}{\rho_t}\r)^{1/2}.\]
Comparing these bounds to \eqref{eq:midtpcII}, we have
\[\PP(D_t<k_t-1) \lesssim  \frac{k_t+ \rho_t}{k_t} \cdot \l(\frac{\rho_t}{m_t}\r)^{1/2} \ \text{ and } \
\PP(D_t>k_t) \lesssim \frac{\rho_t}{k_t} \l(\frac{m_t}{\rho_t}\r)^{1/2}.\]
As a last step, we have to consider the behaviour of $\rho_t/k_t$.
Recall from the previous section that $\PP(D_t\geq \frac 12 m_t)\to 1$ (since $D_t/m_t\to 1$ in probability). 
However, here we have $\PP(D_t\geq k_t) = \PP(N_{k_t}\geq k_t+1)\leq \EE N_{k_t}/(k_t+1) \lesssim (\rho_t/m_t)^{1/2}\to 0$. As a consequence, we must have $k_t> \frac 12 m_t$, and so $\rho_t/k_t < 2 \rho_t/m_t$. This yields that both tail probabilities above are upper bounded by a constant multiple of $(\rho_t/m_t)^{1/2}$. This expression converges to zero, which concludes the proof of the theorem. \hfill $\qed$

\section{The dense regime}\label{section:dense}

As stated in the introduction, \Cref{thm:Gumbel} is a direct consequence of \Cref{thm:Poisson}, hence we only have to concern ourselves with proving the latter. 
We use the notations of the theorem, and in addition, write $A_n^{(k)}$ for the event that $I_n\cup I_{n+1}$ contains an interval of length $r_t$ with at least $n$ points i.e. $A_n^{(k)}=\{X_n^{(k)}=1\}$.

We use a classical result by Arratia, Goldstein and Gordon. 
Define the sets $B_n\coloneq \{n-1,n,n+1\}\cap\{1,\ldots, N\}$: this describes the dependencies between the indicator variables, that is, $X_n^{(k)}$ and $X_m^{(k)}$ are dependent exactly when $m\in B_n$.
Further, define the following quantities:
\begin{align*}
    b_1 & \coloneq \sum_{n=0}^N \sum_{m\in B_n} \PP(A_n^{(k)})\cdot \PP(A_m^{(k)}), \\
    b_2 & \coloneq \sum_{n=0}^N \sum_{n\neq m \in B_n} \PP(A_n^{(k)}\cap A_m^{(k)}), \\
    b_3 & \coloneq \sum_{n=0}^N \EE \l| \EE \l[ X_n-\EE X_n | \{X_m : m\in \{1,\ldots, N\}\setminus B_n\} \r]\r|.
\end{align*}
Then \cite[Theorem~1]{AGG89} states that 
\[\dTV(X^{(k)},Z^{(k)})\leq b_1+b_2+b_3,\]
where, as in the theorem, $Z^{(k)}$ denotes a Poisson random variable with parameter $\EE X^{(k)}=\sum_{n=0}^N \PP(A_i^{(k)})$.

Let $p^{(k)}\coloneq \PP(A_n^{(k)})$ and $q^{(k)}\coloneq \PP(A_{n-1}^{(k)}\cap A_{n}^{(k)})$ for $n=1,\ldots, N$, where clearly the expressions do not depend on the lower index of $A^{(k)}$, with the exception of the boundary case: $\PP(A_n^{(k)})\leq p^{(k)}$ and $\PP(A_{N-1}^{(k)}\cap A_{N}^{(k)})\leq q^{(k)}$. With this notation, we have that $\EE X^{(k)}\sim Np^{(k)}$, as well as the bounds
\[b_1 \leq 3(N+1)(p^{(k)})^2 \ \text{ and } \ b_2 \leq 2(N+1) q^{(k)}.\]
By independence, $b_3=0$. Further, by construction, $N\sim r_t^{-1} = t/\rho_t$.

%We note here that $k$ will be chosen in a way that $\EE X^{(n)} \sim Np^{(n)} \sim p^{(n)}/r_t\to e^{-x}$. That the correct choice is $k_t(x) = a_t +xb_t$ as stated in the theorem will be shown later; however, from this assumption and the computation above, it will follow that $b_1 = O(p^{(n)})$ and $b_2 = O\l(q^{(n)}/p^{(n)}\r).$

The following lemma describes the values of $p^{(k)}$ and $q^{(k)}$ precisely. Note that the expression for $p^{(k)}$ can be found in the book \cite{scan}, originally obtained for the binomial point process in \cite{N65}. Since we use a similar strategy to find $q^{(k)}$, we recall the idea of the proof.
\begin{lemma}\label{lemma:pq}
With the notation above,
\begin{multline*}
    p^{(k)} = 
2\PP(\Po(\rho_t)\geq k)- \PP(\Po(\rho_t)\geq k)^2  \\
    + \PP(\Po(\rho_t)=k) \cdot \l[(k-\rho_t-1) \PP(\Po(\rho_t)\leq k-2) + \rho_t \PP(\Po(\rho_t)= k-2)\r] 
\end{multline*}
and
\[q^{(k)}\leq 3 \PP(\Po(\rho_t)\geq k) +\PP(\Po(\rho_t)=k)^2 \sum_{m=1}^{k-1} \frac{\PP(\Po(\rho_t)\leq m)^2}{\PP(\Po(\rho_t)=m)}.\]
\end{lemma}
\begin{proof}
    First consider $p^{(k)}$ in detail.
    By scaling, we may consider a Poisson point process on $[0,2)$ with intensity $s\coloneq \rho_t$. If $M$ denotes the maximal number of points in an interval of length $1$, then we exactly have $p^{(k)} = \PP(M\geq k)$.
    
    We condition on the number of points in the two halves: let $A_{n,m}$ denote the event where $[0,1]$ and $[1,2]$ contain $n$ and $m$ points, respectively. We may assume that the point $1$ was not chosen, as this happens with probability $0$. 
    Clearly, if $n\geq k$ or $m\geq k$, we have that $\PP(M\geq k \mid A_{n,m})=1$, if $n+m<k$, the probability is zero. In what follows, we assume that $n,m<k$ and $n+m\geq k$.
    Let us now examine $M$ under $A_{n,m}$.
    
    The scanning of the interval $[0,2]$ can be interpreted as follows: we start by considering $[0,1]$, and then sliding this window towards $2$. When the right endpoint of the scanning set hits a point in $(1,2]$, the number of points in it goes up, and when its left endpoint leaves a point in $[0,1]$, it goes down. We can thus consider this scanning as a path of up- and down-steps from $n$ to $m$. This path, and in particular whether it reaches $k$ (corresponding to a $1$-interval with $k$ points), depends only on the order in which the the scanning set leaves and hits the points of the first- and second halves, respectively. What is more, each of these orders occurs with the same probability, and hence the problem becomes purely combinatorial.

    The total number of paths starting at $n$ with $m$ up and $n$ down steps is exactly $\binom{n+m}{n}$; the number of such paths reaching $k$ is $\binom{n+m}{k}$, which can be found using a standard reflection argument that we omit here.
    %each of these paths corresponds to a path with $k$ up and $n+m-k$ down-steps, which can be seen by reflecting the part of the path after the first $k$-hit through the line...

    Putting this together, we have
    \[\PP(M\geq k) = 2\PP(\Po(s)\geq k)- \PP(\Po(s)\geq k)^2 
         + \sum_{\substack{1\leq n,m<k\\ n+m\geq k}} \frac{\binom{n+m}{k}}{\binom{n+m}{n}} \cdot \frac{s^n}{n!} \cdot e^{-s} \cdot \frac{s^m}{m!}\cdot e^{-s},\]
    % \begin{multline*}
    %     \PP(M\geq k) = 2\PP(\Po(s)\geq k)- \PP(\Po(s)\geq k)^2 + \\ 
    %     + \sum_{\substack{1\leq n,m<k\\ n+m\geq k}} \frac{\binom{n+m}{k}}{\binom{n+m}{n}} \cdot \frac{s^n}{n!} \cdot e^{-s} \cdot \frac{s^m}{m!}\cdot e^{-s},
    % \end{multline*}
    where the first part comes from the case $n\geq k$ or $m\geq k$. 
    Examining the sum further, we see by simplifying that the summand is exactly
    \[\frac{e^{-s} s^{k}}{k!} \cdot\frac{e^{-s} s^{n+m-k}}{(n+m-k)!} = \PP(\Po(s) = k) \cdot \PP(\Po(s)=n+m-k),\]
    yielding that after carrying out the substitution $l=n+m-k$, the sum becomes
    \[\PP(\Po(s)=k) \sum_{l=0}^{k-2} \sum_{n=l+1}^{k-1} \PP(\Po(s)=l)= \PP(\Po(s)=k)  \sum_{l=0}^{k-2} (k-l-1) \PP(\Po(s)=l).\]
    By writing $l\PP(\Po(s) = l) = s \PP(\Po(s) = l-1)$, we have
    \begin{multline*}
        \PP(\Po(s)=k)\l[ (k-1) \sum_{l=0}^{k-2} \PP(\Po(s)=l) - s \sum_{l=1}^{k-2} \PP(\Po(s)=l-1)\r] \\
        = \PP(\Po(s)=k) \l[(k-1) \PP(\Po(s)\leq k-2) -s \PP(\Po(s)\leq k-3)\r],
    \end{multline*}
    which simplifies to the expression stated in the lemma.    

    For $q^{(k)}$, we can follow a similar setup, now generating a Poisson point process $\xi$ of intensity $s=\rho_t$ over $[0,3]$. 
    Let $E_1$ denote the event that there is a unit interval in $[0,2]$ with at least $k$ points, and $E_2$ similarly for $[1,3]$. Lastly, let $E=E_1\cap E_2$, thus $\PP(E) = q^{(k)}$. 

    The probability that one of unit intervals $[0,1],[1,2],[2,3]$ contains at least $k$ points is upper bounded by $3\PP(\Po(s)\geq k)$, which contributes the first term in the bound.
    
    Denote by $B_{n,m,n'}$ the event that $\xi$ has $n$, $m$ and $n'$ points in $[0,1]$, $[1,2]$ and $[2,3]$, respectively, where $n,m,n'<k$ and $n+m,n'+m\geq k$.
    We claim that
    \begin{equation}\label{eq:FKGbound}
        \PP(E\mid B_{n,m,n'})\leq \PP(E_1 \mid B_{n,m,n'})\cdot \PP(E_2\mid B_{n,m,n'}) = \frac{\binom{n+m}{k}}{\binom{n+m}{n}} \cdot \frac{\binom{m+n'}{k}}{\binom{m+n'}{m}}.
    \end{equation} 
    We first sketch how this implies the statement of the lemma.
    Since $\PP(B_{n,m,n'}) = \PP(\Po(s)=n)\PP(\Po(s)=m)\PP(\Po(s)=n')$, we have after rearranging that
    \[\PP(E\cap B_{n,m,n'})\leq \PP(\Po(s)=k)^2 \cdot \frac{\PP(\Po(s)=n+m-k) \cdot \PP(\Po(s) = n'+m-k)}{\PP(\Po(s) = m)},\]
    which we must then sum (over $n,m,n'<k$, $n+m,n'+m\geq k$) to obtain the overall probability for $q^{(k)}$.
    With a change of variables for $n+m-k$ and $n'+m-k$, the sum simplifies to
    \[\PP(\Po(s)=k)^2 \sum_{m=1}^{k-1} \frac{\PP(\Po(s)\leq m-1)^2}{\PP(\Po(s)=m)},\]
    completing the proof of the lemma.

    Let us now show \eqref{eq:FKGbound}.
    The equality follows from the combinatorial argument used in the computation of $p^{(k)}$.
    The inequality can be shown with the FKG inequality, which claims the following.
    Let $X_1,\ldots, X_N$ be independent random variables, and $f,g:\R^N\to \R^N$ increasing functions (that is, $f(\mathbf x)\leq f(\mathbf y)$ for all $\mathbf x, \mathbf y\in \R^N$ with $\mathbf x\leq \mathbf y$ componentwise). Then the FKG inequality (see e.g. \cite[Theorem~2.1]{FKG}) states that
    \[\EE[f(X_1,\ldots, X_N)g(X_1,\ldots, X_N)]\geq \EE [f(X_1,\ldots, X_N)] \cdot \EE [g(X_1,\ldots, X_N)]\]
    whenever these expectations exist.
    We use this in our setting the following way. Let $\{X_1,\ldots, X_n\}$, $\{X_{n+1},\ldots, X_{n+m}\}$, $\{X_{n+m+1},\ldots, X_{n+m+n'}\}$ be uniform samples on the intervals $[0,1]$, $[1,2]$, $[2,3]$, respectively, all independent. Set $\mathbf X \coloneq (X_1,\ldots, X_N)$ with $N\coloneq n+m+n'$. %; the support of $\mathbf X$ is then $S\coloneq [0,1]^n \times [1,2]^m \times [2,3]^{n'}$.
    Further, let $f^\ast, g^\ast:\R^n\to \R$ be indicator functions such that $f^\ast(x_1,\ldots, x_N)=1$ whenever there is a unit interval of $[0,2]$ containing at least $k$ of the points $x_1,\ldots, x_N$, and $g^\ast$ similarly for the interval $[1,3]$. 
    With this setup, we have $\PP(E\mid B_{n,m,n'})=\PP(f^*(\mathbf X) = g^*(\mathbf X)=1)$, as well as $\PP(E_1\mid B_{n,m,n'}) =\PP(f^*(\mathbf X) =1)$ and $\PP(E_2\mid B_{n,m,n'}) =\PP(g^*(\mathbf X) =1)$. 
    
    Observe that on the support $S\coloneq [0,1]^n\times [1,2]^m \times [2,3]^{n'}$ of $\bf X$, $f^\ast$ depends only on its first $n+m$ arguments, and in particular is increasing in its first $n$ coordinates and decreasing in the next $m$. Similarly, $g^\ast$ depends only on its last $m+n'$ arguments, is increasing in the middle $m$, and decreasing in the last $n'$ arguments.
    
    Define the transform $\mathbf{h} = (h_1,\ldots, h_N):\R^N\to \R^N$ as
    \[
        h_i(x)\coloneq 
        \begin{cases}
            1-x, & \caseif i=1,\ldots, n\\
            x, & \caseif i=n+1,\ldots, n+m\\
            5-x, & \caseif i=n+m+1,\ldots, N.
        \end{cases}
    \]
    On $S$, applying $\bf h$ gives the reflected point configuration in the intervals $[0,1]$ and $[2,3]$, and the original on $[1,2]$.
    As a consequence, on the support of $\bf X$, the functions $f\coloneq -f^\ast \circ \bf h$ and $g\coloneq g^\ast \circ \bf h$ are increasing in all arguments.
    In addition, since the sample is uniform, $h_i(X_i)$ and $X_i$ have the same law for all $i$, and thus $\mathbf h (\mathbf X) = (h_1(X_1),\ldots, h_N(X_N))$ and $\bf X$ do too.
    Consequently, the pairs $-f(\mathbf X)$ and $f^\ast(\mathbf X)$, and $g(\mathbf X)$ and $g^\ast(\mathbf X)$ also have the same law. %distribution 
    
    Applying the FKG inequality to $f$ and $g$, we obtain the inequality claimed in \eqref{eq:FKGbound}:
    \begin{align*}
        \PP(E \ | \ B_{n,m,n'}) 
        & = \PP(f(\mathbf X)=-1,\ g(\mathbf X)=1) \\
        &\leq  \PP(f(\mathbf X)=1)\PP( g(\mathbf X) = -1) = \PP(E_1 \mid B_{n,m,n'})\cdot \PP(E_2\mid B_{n,m,n'}).
    \end{align*}
\end{proof}

\begin{corollary*}
    With the notation of the proof, 
    \[\PP\l( \frac{M-s}{\sqrt s} \leq x \r)\to (\Phi(x))^2 - (\phi(x))^2 - x\Phi(x)\phi(x)\] as $s\to \infty$ for all $x\in \R$,
    where $\Phi(x)$ and $\phi(x)$ denote the CDF and PDF of the standard normal distribution.
\end{corollary*}
Note that the limiting distribution does not appear to correspond to any classical distribution.
\begin{proof}[Proof of Corollary]
    From \eqref{eq:PoissonAsy}, we have for any $n_s\in \N$ with $n_s=s(1+x/\sqrt s)+o(\sqrt s)$ as $s\to \infty$ that
    \[\PP\l(\Po(s) = n_s\r) \sim \frac{1}{\sqrt{2\pi s}} \cdot e^{-\frac{x^2}{2}} = \frac{\phi(x)}{\sqrt s},\]
    and a standard normal approximation gives that $\PP\l( \Po(s)\leq s\l( 1+\frac{x}{\sqrt s}\r)\r) \to \Phi(x)$. Plugging this into \Cref{lemma:pq} yields the corollary.
\end{proof}

\begin{proof}[Proof of \Cref{thm:Poisson}]

We evaluate the expressions in \Cref{lemma:pq} for $k=k_t\in \N$ as stated in the theorem: let $k_t(x) = a_t+b_t x + o(b_t)$ as $t\to \infty$. 
We will show that the theorem holds for any $(k_t)_{t>0}$ satisfying the asymptotic behaviour, and the precise values are not significant.
Choose $\eps_t(x)$ such that $k_t = \rho_t(1+\eps_t(x))$. Then
\[\eps_t(x) =\sqrt{\frac{2\ln(t/\rho_t)}{\rho_t}} \cdot \l(1+\frac{\ln(\ln(t/\rho_t))}{4\ln (t/\rho_t)} - \frac{\ln \sqrt{\pi}-x}{2\ln(t/\rho_t)} + o\l(\frac{1}{\ln(t/\rho_t)}\r)\r).\]
For most of the computation, we omit the dependence on $x$ to keep the formulas concise.
Observe that $\ln(t/\rho_t)=\ln(r_t^{-d})\to \infty$, and so $\rho_t \eps_t^2\to \infty$, but $\rho_t \eps_t^3\to 0$ as $t\to \infty$.

Start by considering $p^{(k_t)}$.
Using \eqref{eq:PoissonAsy}, we have from \Cref{lemma:pq}
\[\PP(\Po(\rho_t)=k_t)\sim \frac{1}{\sqrt{2\pi \rho_t}}\cdot e^{\rho_t\eps_t} \cdot (1+\eps_t)^{-\rho_t(1+\eps_t)}\sim\frac{1}{\sqrt{2\pi \rho_t}}\cdot \exp \l\{ - \frac{\rho_t\eps_t^2}{2} \r\},\]
where the second expression comes from $\eps_t-(1+\eps_t)\ln(1+\eps_t) = -\eps_t^2/2 + O(\eps_t^3)$: since $\rho_t\eps_t^3\to 0$, the third-order term does not affect the (first-order) asymptotics.

To find the asymptotic behaviour of the tail $\PP(\Po(\rho_t) \geq k_t)$, we may use a result from moderate deviations theory: formally, this can be done by writing a $\Po(\rho_t)$ random variable as the sum of $\lfloor \rho_t \rfloor$ independent Poisson variables with parameter $\rho_t/\lfloor \rho_t \rfloor$, and applying a result of Cramér \cite{C38}; see also the English translation by
Touchette \cite{C38transl}. In particular, by \cite[Theorem~2]{C38transl}, we have
$\PP(\Po(\rho_t)\geq k_t) \sim 1-\Phi(\eps_t \sqrt{\rho_t})$. Note that we again used the $\rho_t\eps_t^3\to 0$ assumption, which can equivalently be written as $(\eps_t \sqrt{\rho_t})^3/\sqrt \rho_t \to 0$, exactly the requirement for the asymptotic equality in \cite{C38transl}.
Since the asymptotic behaviour of $\Phi$ is well known, 
%$1-\Phi(x) \sum e^{-x^2/2} /(\sqrt{2\pi}\x)$ for large $x$
this implies
\begin{equation}\label{eq:epsPoi}
    \PP(\Po(\rho_t) \geq k_t) \sim \frac{1}{\sqrt{2\pi \rho_t \eps_t^2}} \cdot \exp \l\{ - \frac{\rho_t\eps_t^2}{2} \r\}.
\end{equation}
By explicit comparison, and using that $\rho_t\eps_t^2\to \infty$, we see that the third term of the sum in \Cref{lemma:pq} dominates $p^{(k_t)}$, and in particular
\begin{equation}\label{eq:pk}
    p^{(k_t)} \sim (k_t-\rho_t) \PP(\Po(\rho_t)=k_t) \sim \sqrt{\frac{\rho_t \eps_t^2}{2\pi}} \cdot \exp \l\{ - \frac{\rho_t\eps_t^2}{2} \r\}.
\end{equation}
Plugging in the asymptotic expression for $\eps_t$,
\begin{align}\label{eq:eps}
    \frac{\rho_t \eps_t^2(x)}{2} 
    & = \ln \l(\frac{t}{\rho_t}\r) \l(1+\frac{\ln(\ln(t/\rho_t))}{4\ln (t/\rho_t)} - \frac{\ln \sqrt{\pi}-x}{2\ln(t/\rho_t)}  + o\l(\frac{1}{\ln(t/\rho_t)}\r)\r)^2\notag\\
    &  =\ln \l(\frac{t}{\rho_t}\r) \l(1+ \frac{\ln(\ln(t/\rho_t))}{2\ln (t/\rho_t)} - \frac{\ln \sqrt{\pi}-x}{\ln(t/\rho_t)} + O\l(\l(\frac{\ln(\ln(t/\rho_t))}{\ln (t/\rho_t)}\r)^2\r)\r)\notag \\
    & = \ln \l(\frac{t}{\rho_t}\r) + \frac 12 \ln \l(\ln \l(\frac{t}{\rho_t}\r)\r) - \ln \sqrt{\pi} + x + o(1),
\end{align}
where for the negligible terms we use that $\ln(t/\rho_t)\to \infty$.
Plugging this into \eqref{eq:pk}, 
\[p^{(k_t(x))}\sim \sqrt{\frac{\ln(t/\rho_t)}{\pi}} \cdot \l[ \frac{t}{\rho_t} \l(\ln\l(\frac{t}{\rho_t}\r)\r)^{1/2} (\sqrt \pi)^{-1} e^{x}\r]^{-1} \sim \frac{\rho_t}{t} \cdot e^{-x}.\]
Since $\EE X^{(k_t(x)}\sim ( t/\rho_t) p^{(k_t(x))}$, this yields $\EE X^{(k_t(x))}\to e^{-x}$ as $t\to \infty$.
We also immediately have $b_1\lesssim (t/\rho_t) (p^{(k_t)})^2\lesssim \rho_t/t$ in the bounds on the total variation distance.

Lastly, it remains that we bound $q^{(k_t)}$, and thus $b_2\lesssim (t/\rho_t) q^{(k_t)}$. 
From \Cref{lemma:pq}, we have $q^{(k_t)}\leq I_1 + \PP(\Po(\rho_t)=k_t)^2 \cdot (I_2+I_3)$, where $I_1 = 3\PP(\Po(\rho_t)\geq k_t)$ and
\begin{align*}
    %I_1 & = 3\PP(\Po(\rho_t)\geq k_t) \\
    I_2 & = \sum_{m=1}^{K_t} \frac{\PP(\Po(\rho_t)\leq m)^2}{\PP(\Po(\rho_t)=m)}\ \text{ and }\
    I_3 = \sum_{m=K_t+1}^{k_t-1} \frac{\PP(\Po(\rho_t)\leq m)^2}{\PP(\Po(\rho_t)=m)},
\end{align*}
with $K_t\coloneq \lfloor \rho_t\rfloor$ denoting the integer part of $\rho_t$.
The asymptotic of $I_1$ is given in \eqref{eq:epsPoi}.
For $I_2$, we use the bounds \eqref{eq:PoissonUpper} and \eqref{eq:PoissonBounds} on the Poisson probabilities to obtain
\begin{align*}
    I_2 \lesssim \sum_{m=1}^{K_t} e^{-\rho_t}  \sqrt{2\pi m}  \l( \frac{e\rho_t}{m}\r)^m
    &\lesssim\rho_t \sum_{m=1}^{K_t} \frac{e^{-\rho_t}}{\sqrt{2\pi (m-1)}} \l( \frac{e\rho_t}{m-1}\r)^{m-1} \\
    &\lesssim\rho_t \sum_{m=1}^{K_t} \PP(\Po(\rho_t)=m-1) \leq  \rho_t.
\end{align*}
%Alternatively, we could use $\PP(\Po(\rho_t)\leq m)\leq \sqrt m \PP(\Po(\rho_t)=m)$ for $m\leq \rho_t$, which can be derived directly from the Poisson...
Lastly, we turn to the bound for $I_3$, which will be the dominating part of the sum. Upper bounding the nominator by $1$, and using \eqref{eq:PoissonBounds} for the denominator, we have
\[I_3 \lesssim\sum_{m=K_t+1}^{k_t-1} \sqrt{m} \cdot e^{\rho_t}  \cdot  \l(\frac{m}{e\rho_t}\r)^m.\]
The $\sqrt m$ factor is at most some small constant away from $\sqrt \rho_t$ due to the range. Then, by using the monotonicity of the function $y^y (e\rho_t)^{-y}$ for $y\geq \rho_t$, the standard integral bound gives
\[I_3 \lesssim \sqrt{\rho_t} \sum_{m=K_t+1}^{k_t-1} e^{\rho_t} \l(\frac{m}{e\rho_t}\r)^m \lesssim\sqrt \rho_t \int_{\rho_t}^{k_t} e^{\rho_t} \l(\frac{y}{e\rho_t}\r)^y \dd y,\]
where we used that $K_t+1>\rho_t$.
Substituting $y=\rho_t(1+u)$ gives
\[I_3 \lesssim \rho_t^{\frac 32} \int_0^{\eps_t} e^{\rho_t} \l(\frac{1+u}{e}\r)^{\rho_t(1+u)}\dd u = \rho_t^{\frac 32} \int_0^{\eps_t} \l(\frac{\l(1+u\r)^{1+u}}{e^{u} }\r)^{\rho_t}\dd u.\]
Since $(1+u)\ln (1+u)-u = \frac 12 u^2  + O(u^3)$ around $u=0$, and $\eps_t\to 0$ and $\eps^3 \rho_t\to 0$, we have
\[I_3 \lesssim \rho_t^{\frac 32} \cdot \int_0^{\eps_t} e^{\frac{\rho_tu^2}{2}}\cdot e^{O(\rho_t \eps_t^3)}\dd u \lesssim \rho_t^{\frac 32} \cdot \int_0^{\eps_t} e^{\frac{\rho_tu^2}{2}}\dd u.\]
Lastly, let $u=\eps_t \cdot \l(1- \frac{v}{\rho_t\eps_t^2}\r)$ to get
\[I_3 \lesssim \frac{\rho_t^{\frac 12}}{\eps_t} \int_0^{\rho_t\eps_t^2} \exp\l\{\frac{\rho_t\eps_t^2}{2} \l(1-\frac{v}{\rho_t\eps_t^2}\r)^2\r\}\dd v \lesssim \frac{\rho_t^{\frac 12}e^{\frac{\rho_t\eps_t^2}{2}}}{\eps_t} \int_0^{\rho_t\eps_t^2} \exp\l\{-v + \frac{v^2}{2\rho_t\eps_t^2}\r\}\dd v.\]
For $0\leq v\leq \rho_t\eps_t^2$, the bound $-v+v^2/(2\rho_t\eps_t^2)\leq -v/2$ holds, which we may apply to the exponent in the integral. The resulting integral can further be upper bounded by expanding the domain of integration up to infinity, and explicitly evaluating the integral. Overall, this yields
\[I_3 \lesssim \frac{\rho_t^{\frac 12}e^{\frac{\rho_t\eps_t^2}{2}}}{\eps_t} \int_0^{\infty} e^{-\frac v2}\dd v \lesssim \frac{\rho_t^{\frac 12}e^{\frac{\rho_t\eps_t^2}{2}}}{\eps_t} = \rho_t\cdot \frac{ e^{\frac{\rho_t\eps_t^2}{2}}}{\sqrt{\rho_t \eps_t^2}}.\]
Comparing all bounds, we see that the expression involving $I_3$ dominates, and in particular
\[q^{(k_t)} \lesssim \PP(\Po(\rho_t)=k_t)^2 \cdot I_3 \lesssim \frac{e^{-\rho_t\eps_t^2}}{\rho_t} \cdot \frac{\rho_t e^{\frac{\rho_t\eps_t^2}{2}}}{\sqrt{\rho_t \eps_t^2}} = \frac{e^{-\frac{\rho_t\eps_t^2}{2}}}{\sqrt{\rho_t \eps_t^2}}.\]
%Comparing this to \eqref{eq:pk}, $q^{(k_t)}\lesssim p^{(k_t)}/(\rho_t\eps_t^2)$.
Plugging in \eqref{eq:eps}, we have $q^{(k_t)}\lesssim (\rho_t/t)(\ln (t/\rho_t))^{-1}$.
Altogether, the bounds for the total variation distance introduced at the beginning of the section are bounded by $b_1+b_2 \lesssim \l(\ln (t/\rho_t)\r)^{-1}$. Applying \cite[Theorem~1]{AGG89} concludes the proof of the theorem.

\end{proof}

\section*{Funding statement}

This work was funded by the DFG (German Research Foundation) Priority Programme \textit{Random Geometric Systems} under project number 531562368.

\section*{Acknowledgements}

The author thanks Matthias Reitzner for suggesting the topic of this paper, as well as helpful conversations throughout the working process. 
We are also grateful to the anonymous referees for their suggestions that greatly improved the presentation of the paper. We particularly appreciate the referee suggesting an elegant alternative proof to \Cref{lemma:scan}, as well as pointing out an error in the proof of \Cref{lemma:pq}, and proposing the FKG inequality as a solution.

\bibliography{references}

@article{AR20,
  author  = {Akinwande, G. and Reitzner, M.},
  title   = {Multivariate central limit theorems for random simplicial complexes},
  journal = {Adv. Appl. Probab.},
  volume  = {121},
  year    = {2020},
  mrnumber = {4126722},
 note = {\href{https://doi.org/10.1016/j.aam.2020.102076}{DOI:10.1016/j.aam.2020.102076}},
  DOI = {\url{https://doi.org/10.1016/j.aam.2020.102076}},
}

@article{AGG89,
  author  = {Arratia, R. and Goldstein, L. and Gordon, L.},
  title   = {Two moments suffice for Poisson approximations: the Chen--Stein method},
  journal = {Ann. Probab.},
  volume  = {17},
  year    = {1989},
  pages   = {9--25},
  mrnumber = {0972770},
 note = {\href{https://doi.org/10.1214/aop/1176991491}{DOI:10.1214/aop/1176991491}},
  doi = {\url{https://doi.org/10.1214/aop/1176991491}},
}

@article{B16,
  author   = {Bachmann, S.},
  title    = {Concentration for {P}oisson functionals: component counts in random geometric graphs},
  journal  = {Stoch. Process. Appl.},
  volume   = {126},
  year     = {2016},
  pages    = {1306--1330},
  mrnumber = {3473096},
note = {\href{https://doi.org/10.1016/j.spa.2015.11.004}{DOI:10.1016/j.spa.2015.11.004}},
    DOI = {\url{https://doi.org/10.1016/j.spa.2015.11.004}},
}

@article{BP16,
  author   = {Bachmann, S. and Peccati, G.},
  title    = {Concentration bounds for geometric {P}oisson functionals: logarithmic {S}obolev inequalities revisited},
  journal  = {Electron. J. Probab.},
  volume   = {21},
  year     = {2016},
  pages    = {1--44},
  mrnumber = {3485348},
note = {\href{https://doi.org/10.1214/16-EJP4235}{DOI:10.1214/16-EJP4235}},
    DOI = {\url{10.1214/16-EJP4235}},
}

@article{BR17,
  author   = {Bachmann, S. and Reitzner, M.},
  title    = {Concentration for {P}oisson $U$-statistics: subgraph counts in random geometric graphs},
  journal  = {Stoch. Process. Appl.},
  volume   = {128},
  year     = {2017},
  pages    = {3327--3352},
  mrnumber = {3849811},
note = {\href{https://doi.org/10.1016/j.spa.2017.11.001}{DOI:10.1016/j.spa.2017.11.001}},
    DOI={\url{https://doi.org/10.1016/j.spa.2017.11.001}},
}

@article{C38,
  author  = {Cram{\'e}r, H.},
  title   = {Sur un nouveau th{\'e}or{\`e}me limite de la probabilit{\'e}s},
  journal = {Actual. Sci. Indust.},
  volume  = {736},
  year    = {1938},
  pages   = {5--23},
}

@article{C38transl,
  author  = {Cram{\'e}r, H},
  title   = {On a new limit theorem in probability theory},
note = {English translation of ``Sur un nouveau th{\'e}or{\`e}me limite de la probabilit{\'e}s'' by Touchette, H., \href{https://arxiv.org/abs/1802.05988}{arXiv:1802.05988}},
year = {2018},
  URL = {https://arxiv.org/abs/1802.05988},
}

@book{LDPbook,
  author    = {Dembo, A. and Zeitouni, O.},
  title     = {Large Deviations Techniques and Applications},
  edition   = {2},
  publisher = {Springer-Verlag},
  address   = {New York},
  year      = {1998},
  mrnumber  = {1619036},
}

@article{FKG,
    AUTHOR = {Esary, J. D. and Proschan, F. and Walkup, D. W.},
     TITLE = {Association of random variables, with applications},
   JOURNAL = {Ann. Math. Statist.},
  FJOURNAL = {Annals of Mathematical Statistics},
    VOLUME = {38},
      YEAR = {1967},
     PAGES = {1466--1474},
      ISSN = {0003-4851},
   MRCLASS = {60.02 (62.00)},
  MRNUMBER = {217826},
MRREVIEWER = {I.\ Olkin},
note = {\href{https://doi.org/10.1214/aoms/1177698701}{DOI:10.1214/aoms/1177698701}},
       DOI = {10.1214/aoms/1177698701},
       URL = {https://doi.org/10.1214/aoms/1177698701},
}

@article{Gilbert,
  author   = {Gilbert, E. N.},
  title    = {Random plane networks},
  journal  = {J. Soc. Ind. Appl. Math.},
  volume   = {9},
  year     = {1961},
  pages    = {533--543},
}

@book{scan,
  author    = {Glaz, J. and Naus, J. and Wallenstein, S.},
  title     = {Scan Statistics},
  publisher = {Springer},
  address   = {New York},
  year      = {2001},
  mrnumber  = {1869112},
    isbn = {978-1-4419-3167-2},
}

@book{Gruber,
  author    = {Gruber, P. M.},
  title     = {Convex and Discrete Geometry},
  publisher = {Springer},
  address   = {Berlin},
  year      = {2007},
  mrnumber  = {2335496},
    isbn = {978-3-540-71132-2},
}

@article{KTW24,
  author   = {Kahle, M. and Tian, M. and Wang, Y.},
  title    = {On the clique number of noisy random geometric graphs},
  journal  = {Random Struct. Algorithms},
  volume   = {63},
  year     = {2023},
  pages    = {242--279},
  mrnumber = {4605137},
    doi = {\url{https://doi.org/10.1002/rsa.21134}},
note = {\href{https://doi.org/10.1002/rsa.21134}{DOI:10.1002/rsa.21134}},
}

@article{LRP13,
  author   = {Lachi{\`e}ze-Rey, R. and Peccati, G.},
  title    = {Fine Gaussian fluctuations on the {P}oisson space, {I}: contractions, cumulants and geometric random graphs},
  journal  = {Electron. J. Probab.},
  volume   = {18},
  year     = {2013},
  pages    = {1--32},
  mrnumber = {3035760},    
doi = {\url{https://doi.org/10.1214/EJP.v18-2104}},
note = {\href{https://doi.org/10.1214/EJP.v18-2104}{DOI:10.1214/EJP.v18-2104}}}

@article {LM08,
    AUTHOR = {Lao, W. and Mayer, M.},
     TITLE = {{$U$}-max-statistics},
   JOURNAL = {J. Multivariate Anal.},
  FJOURNAL = {Journal of Multivariate Analysis},
    VOLUME = {99},
      YEAR = {2008},
    NUMBER = {9},
     PAGES = {2039--2052},
      ISSN = {0047-259X,1095-7243},
   MRCLASS = {60D05 (60F05 60G70)},
  MRNUMBER = {2466550},
MRREVIEWER = {Pedro\ Ter\'an},
       DOI = {10.1016/j.jmva.2008.02.001},
       URL = {https://doi.org/10.1016/j.jmva.2008.02.001},
note = {\href{https://doi.org/10.1016/j.jmva.2008.02.001}{DOI:10.1016/j.jmva.2008.02.001}}}

@book{PoissonBook,
    AUTHOR = {Last, G\"unter and Penrose, Mathew},
     TITLE = {Lectures on the {P}oisson process},
    SERIES = {Institute of Mathematical Statistics Textbooks},
    VOLUME = {7},
 PUBLISHER = {Cambridge University Press, Cambridge},
      YEAR = {2018},
     PAGES = {xx+293},
      ISBN = {978-1-107-45843-7; 978-1-107-08801-6},
}

@article{M99,
  author   = {M{\aa}nsson, M.},
  title    = {Poisson approximation in connection with clustering of random points},
  journal  = {Ann. Appl. Probab.},
  volume   = {9},
  year     = {1999},
  pages    = {465--492},
  mrnumber = {1687402},
    DOI = {\url{https://doi.org/10.1214/aoap/1029962751}},
note = {\href{https://doi.org/10.1214/aoap/1029962751}{DOI:10.1214/aoap/1029962751}}}

@article{MM11,
  author   = {McDiarmid, C. and M{\"u}ller, T.},
  title    = {On the chromatic number of random geometric graphs},
  journal  = {Combinatorica},
  volume   = {31},
  year     = {2011},
  pages    = {423--488},
  mrnumber = {2861238},
    DOI = {\url{http://dx.doi.org/10.1007/s00493-011-2403-3}},
note = {\href{https://doi.org/10.1007/s00493-011-2403-3}{DOI:10.1007/s00493-011-2403-3}}}

@article{M08,
  author   = {M{\"u}ller, T.},
  title    = {Two-point concentration in random geometric graphs},
  journal  = {Combinatorica},
  volume   = {28},
  year     = {2008},
  pages    = {529--545},
  mrnumber = {2501248},
    DOI={\url{https://doi.org/10.1007/s00493-008-2283-3}},
note = {\href{https://doi.org/10.1007/s00493-008-2283-3}{DOI:10.1007/s00493-008-2283-3}}}

@article{N65,
  author   = {Naus, J.},
  title    = {The distribution of the size of the maximum cluster of points on a line},
  journal  = {J. Amer. Statist. Assoc.},
  volume   = {60},
  year     = {1965},
  pages    = {532--538},
  mrnumber = {0183041},
    doi = {\url{https://doi.org/10.2307/2282688}},
note = {\href{https://doi.org/10.2307/2282688}{DOI:10.2307/2282688}}}

@article{P02,
  author   = {Penrose, M.},
  title    = {Focusing of the scan statistic and geometric clique number},
  journal  = {Adv. Appl. Probab.},
  volume   = {34},
  year     = {2002},
  pages    = {739--753},
  mrnumber = {1938940},
    DOI = {\url{http://dx.doi.org/10.1239/aap/1037990951}},
note = {\href{https://doi.org/10.1239/aap/1037990951}{DOI:10.1239/aap/1037990951}}}

@book{PenroseBook,
  author    = {Penrose, M.},
  title     = {Random Geometric Graphs},
  publisher = {Oxford Univ. Press},
  address   = {Oxford},
  year      = {2003},
  mrnumber  = {1986198},
}

@article{RRvW24,
  author   = {Reitzner, M. and R{\"o}mer, T. and von Westenholz, M.},
  title    = {Covariance matrices of length power functionals of random geometric graphs---an asymptotic analysis},
  journal  = {Linear Algebra Appl.},
  volume   = {691},
  year     = {2024},
  pages    = {151--181},
  mrnumber = {4728617},
    DOI = {\url{https://doi.org/10.1016/j.laa.2024.03.032}},
note = {\href{https://doi.org/10.1016/j.laa.2024.03.032}{DOI:10.1016/j.laa.2024.03.032}}}

@article{RS13,
  author   = {Reitzner, M. and Schulte, M.},
  title    = {Central limit theorems for $U$-statistics of {P}oisson point processes},
  journal  = {Ann. Probab.},
  volume   = {41},
  year     = {2013},
  pages    = {3879--3909},
  mrnumber = {3161465},
    DOI = {\url{https://doi.org/10.1214/12-AOP817}},
note = {\href{https://doi.org/10.1214/12-AOP817}{DOI:10.1214/12-AOP817}}}

@article{RST17,
  author   = {Reitzner, M. and Schulte, M. and Th{\"a}le, C.},
  title    = {Limit theory for the {G}ilbert graph},
  journal  = {Adv. Appl. Probab.},
  volume   = {88},
  year     = {2017},
  pages    = {26--61},
  mrnumber = {3641808},
    DOI = {\url{https://doi.org/10.1016/j.aam.2016.12.006}},
note = {\href{https://doi.org/10.1016/j.aam.2016.12.006}{DOI:10.1016/j.aam.2016.12.006}}}

@article {ST12,
    AUTHOR = {Schulte, Matthias and Th\"ale, Christoph},
     TITLE = {The scaling limit of {P}oisson-driven order statistics with              applications in geometric probability},
   JOURNAL = {Stochastic Process. Appl.},
  FJOURNAL = {Stochastic Processes and their Applications},
    VOLUME = {122},
      YEAR = {2012},
    NUMBER = {12},
     PAGES = {4096--4120},
      ISSN = {0304-4149,1879-209X},
   MRCLASS = {60F17 (60D05 60G55 60H07 62G32)},
  MRNUMBER = {2971726},
       DOI = {10.1016/j.spa.2012.08.011},
       URL = {https://doi.org/10.1016/j.spa.2012.08.011},
note = {\href{https://doi.org/10.1016/j.spa.2012.08.011}{DOI:10.1016/j.spa.2012.08.011}}}

%\begin{bibdiv}
%\biblist{\input{references.ltb}}
%\end{bibdiv}

\end{document}